 \documentclass[preprint,review,12pt]{elsarticle}

 \usepackage{multirow}
 \usepackage{tabularx}
 \usepackage{rotating}

 \usepackage{graphicx}
 \usepackage{epsfig,subfigure,psfrag}

\usepackage{textcomp}
\usepackage{latexsym,amssymb,amsmath}

 \usepackage{lineno}

\begin{document}

\begin{frontmatter}



\title{All-stages-implicit and strong-stability-preserving implicit-explicit Runge-Kutta time discretization schemes for hyperbolic systems with stiff relaxation terms}


\author[a,b]{Shu-Chao Duan\corref{cor1}}\ead{s.duan@163.com}
\author{$et$ $al$}

\cortext[cor1]{Corresponding author}
\address[a]{Institute of Fluid Physics, CAEP, P. O. Box 919-105, Mianyang 621999, China}
\address[b]{Department of Modern Physics, University of Science and Technology of China, Hefei 230026, China}

\begin{abstract}
We construct eight implicit-explicit (IMEX) Runge-Kutta schemes up to third order of the type in which all stages are implicit so that they can be used in the zero relaxation limit in a unified and convenient manner. These all-stages-implicit schemes attain the strong-stability-preserving property in the limiting case, and two are strong-stability-preserving not only for the explicit part but also the implicit part and the entire IMEX scheme. Three schemes can completely recover to the designed accuracy order in two sides of the relaxation parameter for both equilibrium and non-equilibrium initial conditions. Two schemes converge nearly uniformly for equilibrium cases. These all-stages-implicit schemes can be used for hyperbolic systems with stiff relaxation terms or differential equations with some type constraints.
\end{abstract}


\end{frontmatter}



\section{Introduction}\label{sec.introduction}

Several physical phenomena of great importance for applications are described by stiff systems of differential equations of the form $\partial_t U = F(U) + R(U)/\varepsilon$, where $U$ is the vector of conserved quantities, 
$F(U)$ is the vector of hyperbolic fluxes, $R(U)$ denotes the stiff relaxation terms, and $\varepsilon\in[0,\infty]$ is the stiffness parameter, or relaxation parameter in physics. These systems can be related to many types of problem, for example, viscosity/hyper-viscosity, viscoelasticity, heat-conduction, diffusion, turbulence, boundary layer, reacting flows, traffic flows, multiphase flows and phase transitions, kinetic theory of rarefied gases, hydrodynamic models for semiconductors, radiation hydrodynamics, and relaxation magnetohydrodynamics.

One of the major difficulties when computing solutions to the aforementioned problems is the stiffness of the differential equations in temporal integrations. The development of efficient temporal integration schemes for such systems is challenging because the local relaxation time corresponding to the source terms has a wide range and can be much smaller than the global characteristic transport time determined by the characteristic speed and length of the system. In this case, the differential equations are stiff. As a result, if explicit methods are used to integrate the stiff differential equations, the computations become extremely inefficient because the time-step size dictated by the stability requirements is much smaller than that dictated by the characteristic transport speed.

To remove the stability restriction on explicit methods in the case of stiffness, L-stable implicit methods need to be used to filter out the high-frequency component and step over the fast time scale. Pure implicit methods are rarely used because they require the inversion of a large matrix. Implicit-explicit (IMEX) hybrid methods are much more widely used. IMEX methods can be single-step or multi-step. Single-step IMEX methods are also known as additive Runge-Kutta (RK) methods \cite{Ascher:1997, Kennedy:2003}. The fractional step method (or time-splitting method) is another type of hybrid method that can be used. The drawback of using these methods is that it is difficult to exceed second-order temporal accuracy.

Pareschi and Russo \cite{Pareschi:2005} derived IMEX RK schemes up to third order that are strong-stability-preserving (SSP) \cite{Gottlieb:2009} for the limiting system of conservation laws. These schemes, denoted by IMEX-SSP, combine an L-stable implicit and SSP explicit RK scheme into one scheme that satisfies some order conditions. Pareschi and Russo's IMEX-SSP schemes are widely used in studies because they have the asymptotic preserving property, that is, the consistency of the scheme with the equilibrium system and asymptotic accuracy, thus, the order of accuracy is maintained in the stiff limit \cite{Pareschi:2005}.

However, these schemes have a minor defect that makes them cumbersome when used in the zero relaxation limit because the relaxation parameter $\varepsilon$ is a denominator in the final explicit assembly stage of a standard RK scheme. The denominator definitely should not be zero. Even when the denominator is not zero but finite small, the machine truncation error would be large. Therefore, we need an IMEX scheme whose final stage is also implicit or that has no final assembly stage so that the zero or small parameter can be multiplied on the left-hand side. This is possible, in fact, the second type of Ascher-Ruuth-Spiteri (ARS) \cite{Ascher:1997} IMEX scheme is of this type. Liu and Zou \cite{Liu:2006} thoroughly studied this type of IMEX scheme, motivated by the consideration of the convenience of enforcing some type constraints, such as divergence-free. They provided such schemes up to fourth order. We note that Liu and Zou's schemes have some differences compared with the second type of ARS.

The aim of the present paper is to construct new IMEX RK schemes by combining an L-stable implicit and SSP explicit RK scheme, and ensuring that the final assembly stage of the implicit part is also implicit or is not required. Thus, all stages are implicit so that the scheme can be used in the zero relaxation limit in a unified and convenient manner. All the properties of Pareschi and Russo's IMEX-SSP schemes are held. Optimizations are performed to maximize the absolute monotonicity region \cite{Higueras:2006, Higueras:2014}, the intersection of the stable region and the imaginary axis, or stable region of the IMEX RK schemes.

\section{IMEX RK schemes}

An IMEX RK scheme consists of applying an implicit RK scheme to the stiff source terms and an explicit scheme to the nonstiff terms. When applied to $\partial_t U = R(U)/\varepsilon+F(U)$, an $s$-stage IMEX RK scheme takes the following form \cite{Kennedy:2003, Pareschi:2005}:
\begin{eqnarray}
U^{(i)} &=& U^n + \Delta t \sum_{j=1}^s a_{ij}\frac{R(U^{(j)})}{\varepsilon} 
          + \Delta t \sum_{j=1}^{i-1} b_{ij}F(U^{(j)}),\label{eq.Ui}\\
U^{(n+1)} &=& U^n + \Delta t \sum_{i=1}^s w_i\frac{R(U^{(i)})}{\varepsilon} 
          + \Delta t \sum_{i=1}^{s} \omega_iF(U^{(i)}).\label{eq.Un1}
\end{eqnarray}
The matrices $A=(a_{ij})$ and $B=(b_{ij})$, where $b_{ij}=0$ for $j\ge i$, are coefficient matrices and correspond to the implicit and explicit part of the entire mixed scheme, respectively. The vectors $w = (w_1, ..., w_s)^T$ and $\omega=(\omega_1,...,\omega_s)^T$ are the weight vectors. Note that we placed the implicit and explicit parts in a different sequence compared with some papers because we have written them in the order indicated by the term IMEX itself: IM first, then EX. Typically, we only consider diagonally implicit RK (DIRK) schemes \cite{Kennedy:2003} for the implicit part ($a_{ij}=0$ for $j>i$) because of their simplicity and efficiency when solving algebraic equations. As standard, we can represent the IMEX RK scheme using a Butcher double tableau:
\begin{equation}
\centering
\begin{tabular}{c|c}
$c$ & $A$ \\
\hline
  & $w$
\end{tabular}
\ \ 
\begin{tabular}{c|c}
$d$ & $B$ \\
\hline
  & $\omega$
\end{tabular}
\end{equation}
where vectors $c=(c_1, ..., c_s)^T$ and $d=(d_1, ..., d_s)^T$ specify the internal sub-time level.

We refer to the stages corresponding to Eq. (\ref{eq.Ui}) as the intermediate stages and the stage corresponding to Eq. (\ref{eq.Un1}) as the final assembly stage. The final assembly stage is explicit for both the explicit and implicit part, although the scheme is an IMEX RK scheme. As noted previously, this type of scheme has the aforementioned minor defect. The scheme is stiffly accurate (SA) if $w_i = a_{si}$. If we let $w_i = a_{si}$ and $\omega_i = b_{si}$, then the final assembly stage is not needed because $U^{n+1} = U^{(s)}$ and the Butcher double tableau can be simplified to the following, in which weight vectors are no longer necessary \cite{Ascher:1997, Liu:2006}:
\begin{equation}\label{eq.ASIStyle}
\centering
\begin{tabular}{c|c}
$c$ & $A$ \\
\hline
\end{tabular}
\ \ 
\begin{tabular}{c|c}
$d$ & $B$ \\
\hline
\end{tabular}
\end{equation}

In this case, all stages are implicit (all-stages-implicit (ASI) type) so that the zero or small parameter can be multiplied on the left-hand side of Eq. (\ref{eq.Ui}). Thus, the resulting scheme could be used in the zero relaxation limit in a unified and convenient manner, for example, the second type of ARS \cite{Ascher:1997} IMEX scheme and Liu and Zou's scheme \cite{Liu:2006}. However, the explicit parts of both these schemes are not SSP and the implicit parts of Liu and Zou's scheme are explicit singly diagonally implicit RK (ESDIRK) \cite{Kennedy:2003}, characterized by having an explicit first stage. Thus, these schemes do not have the good properties of Pareschi and Russo's schemes \cite{Pareschi:2005} for the application to hyperbolic systems with stiff relaxation terms.

\section{New schemes}

We construct new ASI type (i.e., form (\ref{eq.ASIStyle})) IMEX RK schemes up to third order. First, we fix the explicit part by a known optimal \cite{Gottlieb:2009} or optimized \cite{Higueras:2014} SSP explicit RK scheme and then impose order conditions described in studies such as those conducted by Pareschi and Russo \cite{Pareschi:2005} and Liu and Zou \cite{Liu:2006}. Some degrees of freedom (DOFs) remain. The implicit part that we considered is SA DIRK or zeroed in the first column, so it has a vanishing stability function at minus infinity, that is, $R(-\infty)=0$, which makes it L-stable if it is A-stable, according to Proposition 3.8 in the study by Hairer and Wanner \cite{Hairer:1987}. Optimizations are performed in the remaining DOF space to maximize the absolute monotonicity region \cite{Higueras:2006, Higueras:2014}, the intersection of the stable region and the imaginary axis, or the area of the stable region of the IMEX RK schemes, and simultaneously, the remaining DOFs are fixed. The schemes presented below are named ASI-SSP($m$,$n$,$p$), where $m$ denotes the number of stages of the implicit part, $n$ the explicit part, and $p$ the accuracy order of the IMEX schemes. The first is ASI-SSP(4,3,2):
\begin{equation}
\centering
\begin{tabular}{c|cccc}
$c_1$ & $\gamma$ & 0 & 0 & 0 \\
$c_2$ & $\alpha$ & $\gamma$ & 0 & 0 \\
$c_3$ & $\beta$ & $a_{32}$ & $\gamma$ & 0 \\
$c_4$ & $a_{41}$ & $a_{42}$ & $a_{43}$ & $\gamma$ \\
\hline
\end{tabular}
\ \ 
\begin{tabular}{c|cccc}
0 & 0 & 0 & 0 & 0 \\
$\frac{1}{2}$ & $\frac{1}{2}$ & 0 & 0 & 0 \\
1 & $\frac{1}{2}$ & $\frac{1}{2}$ & 0 & 0 \\
1 & $\frac{1}{3}$ & $\frac{1}{3}$ & $\frac{1}{3}$ & 0 \\
\hline
\end{tabular}
\end{equation}
where
\begin{equation}
\begin{split}
a_{32} &= \frac{3}{2} - \alpha - \beta - 3\gamma, \\
a_{41} &= \frac{-1+\gamma+2(1+\gamma)\alpha + 4\gamma^2}{3(-1+2\alpha+2\gamma)}, \\
a_{42} &= 1 - 2a_{41}, \\
a_{43} &= a_{41} - \gamma.
\end{split}
\end{equation}
We fix $\gamma=1/4$ to provide L-stability and a small error constant for the above implicit RK scheme \cite{Hairer:1987}. Then, the remaining two DOFs are fixed by maximizing the region of absolute monotonicity $R(A,B)$ of the IMEX scheme, which is defined in the studies by Higueras \cite{Higueras:2006} and
Higueras et al. \cite{Higueras:2014}. A detailed scan shows that the largest value of $r_1$ such that $(r_1,0)\in R(A,B)$ is $r_1 = 2(\sqrt{5}-1)$. This value is attained for $1/2 \le \alpha \le (3+\sqrt{5})/8$ and $\beta=3/4-\alpha$. One DOF still remains. We then maximize the area of the stable region of the IMEX scheme. The definition of the stable region of an IMEX scheme and the method to determine it are provided in reference \cite{Duan:2016}. We find that the maximum area is approximately 14.57 when $\alpha=1/2$. As a comparison, the area of the stable region of the explicit part SSP(3,2) is approximately 16.05. The completely determined ASI-SSP(4,3,2) is
\begin{equation}
\centering
\begin{tabular}{c|cccc}
$\frac{1}{4}$ & $\frac{1}{4}$ & 0             & 0 & 0 \\
$\frac{3}{4}$ & $\frac{1}{2}$ & $\frac{1}{4}$ & 0 & 0 \\
$\frac{1}{2}$ & $\frac{1}{4}$ & 0             & $\frac{1}{4}$ & 0 \\
1             & $\frac{1}{2}$ & 0             & $\frac{1}{4}$ & $\frac{1}{4}$ \\
\hline
\end{tabular}
\ \ 
\begin{tabular}{c|cccc}
0 & 0 & 0 & 0 & 0 \\
$\frac{1}{2}$ & $\frac{1}{2}$ & 0 & 0 & 0 \\
1 & $\frac{1}{2}$ & $\frac{1}{2}$ & 0 & 0 \\
1 & $\frac{1}{3}$ & $\frac{1}{3}$ & $\frac{1}{3}$ & 0 \\
\hline
\end{tabular}
\end{equation}

The explicit part of ASI-SSP(4,3,2) is the known optimal SSP(3,2). The implicit part and the entire IMEX scheme are also SSP.

We note that there are no (2,2,2) and (3,2,2) combinations with non-zero diagonal entries. There is a (4,2,2) combination that fulfils all requisite conditions, but it is not competitive with ASI-SSP(4,3,2) as a four-stage scheme because of the smaller area of the stable region, and hence has not been presented in this paper.

If we let the first column of the (4,3,2) combination be zeros, we obtain a three-stage scheme ASI-SSP(3',3,2). The prime for the first number 3 indicates that the three-stage implicit RK is used in a nominal four-stage form. ASI-SSP(3',3,2) has two DOFs (vector $d = c$ is no longer necessary):
\begin{equation}
\centering
\begin{tabular}{c|cccc|cccc}
0 & 0 & 0 & 0 & 0 & 0 & 0 & 0 & 0 \\
$\frac{1}{2}$ & 0 & $\frac{1}{2}$ & 0 & 0 & $\frac{1}{2}$ & 0 & 0 & 0 \\
1 & 0 & $1-\alpha$  & $\alpha$ & 0 & $\frac{1}{2}$ & $\frac{1}{2}$ & 0 & 0 \\
1 & 0 & 1 & $-\beta$ & $\beta$ & $\frac{1}{3}$ & $\frac{1}{3}$ & $\frac{1}{3}$ & 0 \\
\hline
\end{tabular}
\end{equation}

This scheme is not SSP because of the negative entry \cite{Higueras:2006, Higueras:2014}, but it attains an SSP property in the same manner as Pareschi and Russo's IMEX-SSP schemes \cite{Pareschi:2005}. If we impose a singly diagonal condition, the resulting SDIRK would be more efficient when Newton-type iterative methods are used to obtain a solution. This ASI-SSP(3',3,2) is
\begin{equation}
\centering
\begin{tabular}{c|cccc|cccc}
0 & 0 & 0 & 0 & 0 & 0 & 0 & 0 & 0 \\
$\frac{1}{2}$ & 0 & $\frac{1}{2}$ & 0 & 0 & $\frac{1}{2}$ & 0 & 0 & 0 \\
1 & 0 & $\frac{1}{2}$ & $\frac{1}{2}$ & 0 & $\frac{1}{2}$ & $\frac{1}{2}$ & 0 & 0 \\
1 & 0 & 1 & $\frac{-1}{2}$ & $\frac{1}{2}$ & $\frac{1}{3}$ & $\frac{1}{3}$ & $\frac{1}{3}$ & 0 \\
\hline
\end{tabular}
\end{equation}
The area of its stable region is approximately 11.54. If we let $\alpha=2/25, \beta=3/8$, we obtain another ASI-SSP(3',3,2) with a larger area for the stable region of approximately 12.80, which is very close to the maximum:
\begin{equation}
\centering
\begin{tabular}{c|cccc|cccc}
0 & 0 & 0 & 0 & 0 & 0 & 0 & 0 & 0 \\
$\frac{1}{2}$ & 0 & $\frac{1}{2}$ & 0 & 0 & $\frac{1}{2}$ & 0 & 0 & 0 \\
1 & 0 & $\frac{23}{25}$ & $\frac{2}{25}$ & 0 & $\frac{1}{2}$ & $\frac{1}{2}$ & 0 & 0 \\
1 & 0 & 1 & $\frac{-3}{8}$ & $\frac{3}{8}$ & $\frac{1}{3}$ & $\frac{1}{3}$ & $\frac{1}{3}$ & 0 \\
\hline
\end{tabular}
\end{equation}

We substitute the optimal SSP(3,2) in ASI-SSP(4,3,2) with an optimized SSP(3',2) \cite{Higueras:2014} and we obtain ASI-SSP(4,3',2). The prime for the number 3 discriminates the optimized SSP(3',2) from the optimal SSP(3,2). ASI-SSP(4,3',2) has four DOFs:
\begin{equation}
\centering
\begin{tabular}{c|cccc}
$c_1$ & $\gamma$ & 0 & 0 & 0 \\
$c_2$ & $a_{12}$ & $\gamma$ & 0 & 0 \\
$c_3$ & $\alpha$ & $a_{32}$ & $\gamma$ & 0 \\
$c_4$ & $a_{41}$ & $a_{42}$ & $\beta$ & $\gamma$ \\
\hline
\end{tabular}
\ \ 
\begin{tabular}{c|cccc}
0 & 0 & 0 & 0 & 0 \\
$\frac{5}{6}$ & $\frac{5}{6}$ & 0 & 0 & 0 \\
$\frac{1}{3\delta}$ & $\frac{1}{6\delta}$ & $\frac{1}{6\delta}$ & 0 & 0 \\
1 & $\frac{4}{5}-\delta$ & $\frac{1}{5}$ & $\delta$ & 0 \\
\hline
\end{tabular}
\end{equation}
where
\begin{equation}
\begin{split}
a_{21} &= \frac{5(\beta/2-\delta/8)}{6\beta-3\delta}, \\
a_{32} &= \frac{5/2-2a_{21}-10\alpha\delta}{10\delta}, \\
a_{41} &= \frac{9}{20} - \beta + \frac{2\beta}{5\delta}, \\
a_{42} &= \frac{3}{10} - \frac{2\beta}{5\delta}.
\end{split}
\end{equation}
As previously, we fix $\gamma = 1/4$, and then three DOFs remain. We find that $\alpha=(391-36\sqrt{5})/840, \beta=0, \delta = 7/25$ provides the maximum $r_1 = 2(\sqrt{5}-1)$ and an area of the stable region of approximately 10.70, which is very close to the maximum, and is the area of the stable region of the explicit part SSP(3',2). The implicit part and the entire IMEX scheme are both SSP. Using these values, the ASI-SSP(4,3',2) becomes
\begin{equation}
\centering
\begin{tabular}{c|cccc}
$\frac{1}{4}$ & $\frac{1}{4}$ & 0 & 0 & 0 \\
$\frac{11}{24}$ & $\frac{5}{24}$ & $\frac{1}{4}$ & 0 & 0 \\
$\frac{167}{168}$ & $\frac{391-36\sqrt{5}}{840}$ & $\frac{3(13+2\sqrt{5})}{140}$ & $\frac{1}{4}$ & 0 \\
1 & $\frac{9}{20}$ & $\frac{3}{10}$ & 0 & $\frac{1}{4}$ \\
\hline
\end{tabular}
\ \ 
\begin{tabular}{c|cccc}
0 & 0 & 0 & 0 & 0 \\
$\frac{5}{6}$ & $\frac{5}{6}$ & 0 & 0 & 0 \\
$\frac{25}{21}$ & $\frac{25}{42}$ & $\frac{25}{42}$ & 0 & 0 \\
1 & $\frac{13}{25}$ & $\frac{1}{5}$ & $\frac{7}{25}$ & 0 \\
\hline
\end{tabular}
\end{equation}

We let the first column of the (4,3',2) combination be zeros. We obtain a three-stage scheme ASI-SSP(3',3',2) with one DOF:
\begin{equation}
\centering
\begin{tabular}{c|cccc|cccc}
0 & 0 & 0 & 0 & 0 & 0 & 0 & 0 & 0 \\
$\frac{5}{6}$ & 0 & $\frac{5}{6}$ & 0 & 0 & $\frac{5}{6}$ & 0 & 0 & 0 \\
$\frac{1}{3\delta}$ & 0 & $\frac{1}{3\delta}-\frac{5}{6}$ & $\frac{5}{6}$ & 0 & $\frac{1}{6\delta}$ & $\frac{1}{6\delta}$ & 0 & 0 \\
1 & 0 & $\frac{1+6\delta}{3(2-5\delta)}$ & $\frac{-17\delta}{6(2-5\delta)}$ & $\frac{5}{6}$ & $\frac{4}{5}-\delta$ & $\frac{1}{5}$ & $\delta$ & 0 \\
\hline
\end{tabular}
\end{equation}
This scheme cannot be SSP because all entries cannot be non-negative simultaneously. When $\delta=1/5$, the area of the stable region achieves a maximum of approximately 8.77. Using this value, the ASI-SSP(3',3',2) becomes
\begin{equation}
\centering
\begin{tabular}{c|cccc|cccc}
0 & 0 & 0 & 0 & 0 & 0 & 0 & 0 & 0 \\
$\frac{5}{6}$ & 0 & $\frac{5}{6}$ & 0 & 0 & $\frac{5}{6}$ & 0 & 0 & 0 \\
$\frac{5}{3}$ & 0 & $\frac{5}{6}$ & $\frac{5}{6}$ & 0 & $\frac{5}{6}$ & $\frac{5}{6}$ & 0 & 0 \\
1 & 0 & $\frac{11}{15}$ & $\frac{-17}{30}$ & $\frac{5}{6}$ & $\frac{5}{3}$ & $\frac{1}{5}$ & $\frac{1}{5}$ & 0 \\
\hline
\end{tabular}
\end{equation}

Another solution of the (3',3',2) combination with a smaller maximum area of the stable region is omitted here.

The final second-order scheme we present is ASI-SSP(4',4,2):
\begin{equation}
\centering
\begin{tabular}{c|ccccc|ccccc}
0 & 0 & 0 & 0 & 0 & 0 & 0 & 0 & 0 & 0 & 0\\
$\frac{1}{3}$ & 0 & $\frac{1}{3}$ & 0 & 0 & 0 & $\frac{1}{3}$ & 0 & 0 & 0 & 0\\
$\frac{2}{3}$ & 0 & $\frac{1}{3}$ & $\frac{1}{3}$ & 0 & 0 & $\frac{1}{3}$ & $\frac{1}{3}$ & 0 & 0 & 0\\
1 & 0 & $\alpha$ & $\frac{2}{3} - \alpha$ & $\frac{1}{3}$ & 0 & $\frac{1}{3}$ & $\frac{1}{3}$ & $\frac{1}{3}$ & 0 & 0\\
1 & 0 & $\beta$ & $\frac{3}{2} - 2\beta$ & $\beta-\frac{5}{6}$ & $\frac{1}{3}$ & $\frac{1}{4}$ & $\frac{1}{4}$ & $\frac{1}{4}$ & $\frac{1}{4}$ & 0\\
\hline
\end{tabular}
\end{equation}
This scheme cannot be SSP because all entries cannot be non-negative simultaneously. When $\alpha=1/5, \beta=1/2$, the area of the stable region achieves a local maximum of approximately 27.84. As a comparison, the area of the stable region of the explicit part SSP(4,2) is approximately 32.26. Using this value, the scheme becomes
\begin{equation}
\centering
\begin{tabular}{c|ccccc|ccccc}
0 & 0 & 0 & 0 & 0 & 0 & 0 & 0 & 0 & 0 & 0\\
$\frac{1}{3}$ & 0 & $\frac{1}{3}$ & 0 & 0 & 0 & $\frac{1}{3}$ & 0 & 0 & 0 & 0\\
$\frac{2}{3}$ & 0 & $\frac{1}{3}$ & $\frac{1}{3}$ & 0 & 0 & $\frac{1}{3}$ & $\frac{1}{3}$ & 0 & 0 & 0\\
1 & 0 & $\frac{1}{5}$ & $\frac{7}{15}$ & $\frac{1}{3}$ & 0 & $\frac{1}{3}$ & $\frac{1}{3}$ & $\frac{1}{3}$ & 0 & 0\\
1 & 0 & $\frac{1}{2}$ & $\frac{1}{2}$ & $\frac{-1}{3}$ & $\frac{1}{3}$ & $\frac{1}{4}$ & $\frac{1}{4}$ & $\frac{1}{4}$ & $\frac{1}{4}$ & 0\\
\hline
\end{tabular}
\end{equation}
When $\alpha=10/9, \beta=6/5$, the area of the stable region achieves a local maximum of approximately 27.86. Using this value, the scheme becomes
\begin{equation}
\centering
\begin{tabular}{c|ccccc|ccccc}
0 & 0 & 0 & 0 & 0 & 0 & 0 & 0 & 0 & 0 & 0\\
$\frac{1}{3}$ & 0 & $\frac{1}{3}$ & 0 & 0 & 0 & $\frac{1}{3}$ & 0 & 0 & 0 & 0\\
$\frac{2}{3}$ & 0 & $\frac{1}{3}$ & $\frac{1}{3}$ & 0 & 0 & $\frac{1}{3}$ & $\frac{1}{3}$ & 0 & 0 & 0\\
1 & 0 & $\frac{10}{9}$ & $\frac{-4}{9}$ & $\frac{1}{3}$ & 0 & $\frac{1}{3}$ & $\frac{1}{3}$ & $\frac{1}{3}$ & 0 & 0\\
1 & 0 & $\frac{6}{5}$ & $\frac{-9}{10}$ & $\frac{11}{30}$ & $\frac{1}{3}$ & $\frac{1}{4}$ & $\frac{1}{4}$ & $\frac{1}{4}$ & $\frac{1}{4}$ & 0\\
\hline
\end{tabular}
\end{equation}

The first third-order scheme is ASI-SSP(6,4,3): 
\begin{equation}
\centering
\begin{tabular}{c|cccccc|cccccc}
$\frac{1}{3}$ & $\frac{1}{3}$ & 0 & 0 & 0 & 0 & 0 & 0 & 0 & 0 & 0 & 0 & 0\\
0 & $\frac{-1}{3}$ & $\frac{1}{3}$ & 0 & 0 & 0 & 0 & 0 & 0 & 0 & 0 & 0 & 0\\
$\frac{1}{2}$ & $\frac{1}{6}-\alpha$ & $\alpha$ & $\frac{1}{3}$ & 0 & 0 & 0 & 0 & $\frac{1}{2}$ & 0 & 0 & 0 & 0\\
1 & $\frac{1}{6}-2\alpha$ & $2\alpha$ & $\frac{1}{2}$ & $\frac{1}{3}$ & 0 & 0 & 0 & $\frac{1}{2}$ & $\frac{1}{2}$ & 0 & 0 & 0\\
$\frac{1}{2}$ & $\alpha$ & $\frac{1}{3}-\alpha+\beta$ & $\frac{-1}{6}-2\beta$ & $\beta$ & $\frac{1}{3}$ & 0 & 0 & $\frac{1}{6}$ & $\frac{1}{6}$ & $\frac{1}{6}$ & 0 & 0\\
1 & 0 & $\frac{1}{6}$ & $\frac{1}{2}$ & $\frac{-1}{6}$ & $\frac{1}{6}$ & $\frac{1}{3}$ & 0 & $\frac{1}{6}$ & $\frac{1}{6}$ & $\frac{1}{6}$ & $\frac{1}{2}$ & 0\\
\hline
\end{tabular}
\end{equation}
The singly diagonal entries have been selected to minimize the error constant \cite{Hairer:1987}. We maximize the intersection of the stable region and the imaginary axis. When $\alpha=-3/10, \beta=-7/10$, the intersection approaches the maximum, but the area of the stable region, approximately 3.98, is away from the maximum. If we maximize the area, when $\alpha=14/25, \beta=-3/25$, the area of the stable region is approximately 18.34, which is very close to the maximum. As a comparison, the area of the stable region of the explicit part SSP(4,3) is approximately 19.61.

The second third-order scheme is ASI-SSP(5',4,3):
\begin{equation}
\centering
\begin{tabular}{c|cccccc|cccccc}
0 &0 & 0 & 0 & 0 & 0 & 0 & 0 & 0 & 0 & 0 & 0 & 0\\
$\frac{1}{2}$ & 0 & $\frac{1}{2}$ & 0 & 0 & 0 & 0 & $\frac{1}{2}$ & 0 & 0 & 0 & 0 & 0\\
1 & 0 & $\frac{1}{2}$ & $\frac{1}{2}$ & 0 & 0 & 0 & $\frac{1}{2}$ & $\frac{1}{2}$ & 0 & 0 & 0 & 0\\
$\frac{1}{2}$ & 0 & $\frac{1}{2}$ & $\frac{-1}{2}$ & $\frac{1}{2}$ & 0 & 0 & $\frac{1}{6}$ & $\frac{1}{6}$ & $\frac{1}{6}$ & 0 & 0 & 0\\
0 & 0 & $-1-\alpha$ & $\frac{1}{2}$ & $\alpha$ & $\frac{1}{2}$ & 0 & 0 & 0 & 0 & 0 & 0 & 0\\
1 & 0 & $\frac{2}{3}$ & $\frac{-1}{3}$ & 0 & $\frac{1}{6}$ & $\frac{1}{2}$ & $\frac{1}{6}$ & $\frac{1}{6}$ & $\frac{1}{6}$ & $\frac{1}{2}$ & 0 & 0\\
\hline
\end{tabular}
\end{equation}
The maximum area of the stable region of approximately 14.22 is attained for $\alpha=-3$. For this scheme, the optimal SSP(4,3) is zero-padded in the fifth row and column to create a five-stage form. There are two other padding methods: in the third or fourth row and column, that provide smaller maximum areas for the stable region, hence they are omitted here.

The final third-order scheme is ASI-SSP(5',5,3): 

\begin{equation}
\centering
\newcommand{\mycom}[1]{{\tiny\begin{tabular}{@{}r@{}}#1\end{tabular}}}
\begin{sideways}
\begin{tabular}{r|rrrrrr|rrrrrr}
0 & 0 & 0 & 0 & 0 & 0 & 0 & 0 & 0 & 0 & 0 & 0 & 0\\
\mycom{0.37726891\\53313681} & 0 & \mycom{0.37726891\\53313681} & 0 & 0 & 0 & 0 & \mycom{0.37726891\\53313681} & 0 & 0 & 0 & 0 & 0\\
\mycom{0.75453783\\06627362} & 0 & \mycom{0.37726891\\53313681} & \mycom{0.37726891\\53313681} & 0 & 0 & 0 & \mycom{0.37726891\\53313681} & \mycom{0.37726891\\53313681} & 0 & 0 & 0 & 0\\
\mycom{0.72898566\\16121875} & 0 & \mycom{0.75280719\\94958022} & \mycom{-0.40109045\\321498277} & \mycom{0.37726891\\53313681} & 0 & 0 & \mycom{0.24299522\\053739583} & \mycom{0.24299522\\053739583} & \mycom{0.24299522\\053739583} & 0 & 0 & 0\\
\mycom{0.69922613\\59316696} & 0 & \mycom{0.98566914\\07902044} & \mycom{-0.41371192\\01899029} & -0.25 & \mycom{0.37726891\\53313681} & 0 & \mycom{0.15358906\\769512654} & \mycom{0.15358906\\769512654} & \mycom{0.15358906\\769512654} & \mycom{0.23845893\\284629002} & 0 & 0\\
1 & 0 & \mycom{0.86304437\\29722925} & \mycom{0.33350418\\45496738} & \mycom{-1.79042099\\09531452} & \mycom{1.21660351\\8099811} & \mycom{0.37726891\\53313681} & \mycom{0.20673402\\086480455} & \mycom{0.20673402\\086480455} & \mycom{0.11709725\\184184275} & \mycom{0.18180256\\01201412} & \mycom{0.28763214\\630840694} & 0\\
\hline
\end{tabular}
\end{sideways}
\end{equation}
The area of the stable region is approximately 17.96, whereas it is approximately 33.49 for the optimal SSP(5,3) scheme.

\begin{figure}[htbp]
\centering
  \includegraphics[width=0.4\textwidth]{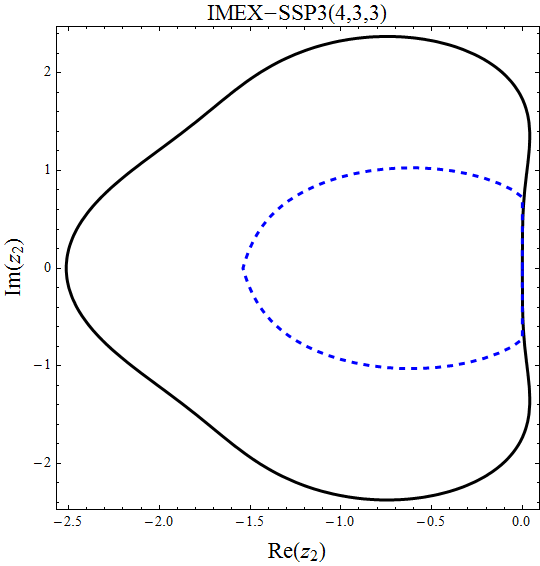}
  \includegraphics[width=0.4\textwidth]{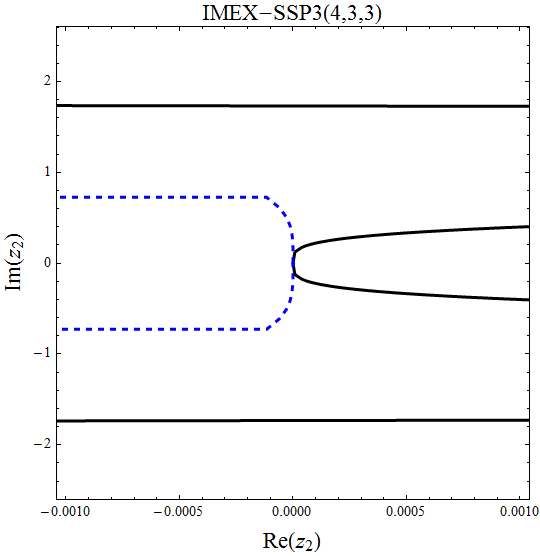}\\
\caption{Stable regions for eight schemes designed in this paper and two from Pareschi and Russo [3]. The blue dashed lines are for IMEX schemes and the black solid lines are for the explicit part. The right column contains close-ups of images on the left.}\label{fig.StableRegions}
\end{figure}

According to Higueras et al. \cite{Higueras:2014}, for methods aimed at hyperbolic systems, it is better to make the stable region contain an interval on the imaginary axis or at least be sufficiently large close to the imaginary axis. All second-order optimal SSP explicit RK schemes have a stable region that contains no part of the imaginary axis, which is also the case for the IMEX schemes based on such explicit RK schemes. Figure \ref{fig.StableRegions} shows the stable regions for eight schemes designed in this paper and two from Pareschi and Russo \cite{Pareschi:2005} for comparison. Only two of these 10 schemes, ASI-SSP(4,3',2) and ASI-SSP(6,4,3), can be optimized to have a nontrivial interval on the imaginary axis. ASI-SSP(4,3',2) is noteworthy because there seemingly exists a point to maximize $r_1$, the area of the stable region, and the interval on the imaginary axis, simultaneously, and $\alpha=(391-36\sqrt{5})/840, \beta=0, \delta=7/25$ is very close to this point. ASI-SSP(6,4,3) with $\alpha=-3/10, \beta=-7/10$ has an interval of slightly more than 1.1 on the imaginary axis. Some of the remaining schemes seemingly approach the imaginary axis closely, but do not contain any of it. 

The precedence in optimization is first, absolute monotonicity, then the intersection, and finally, the area. We cannot confirm that this is the best strategy. Thus, we present the undefined forms with free parameters together with the defined forms in the manuscript so that they may be re-optimized in future research.

\begin{figure}[htbp]
\centering
  \includegraphics[width=0.4\textwidth]{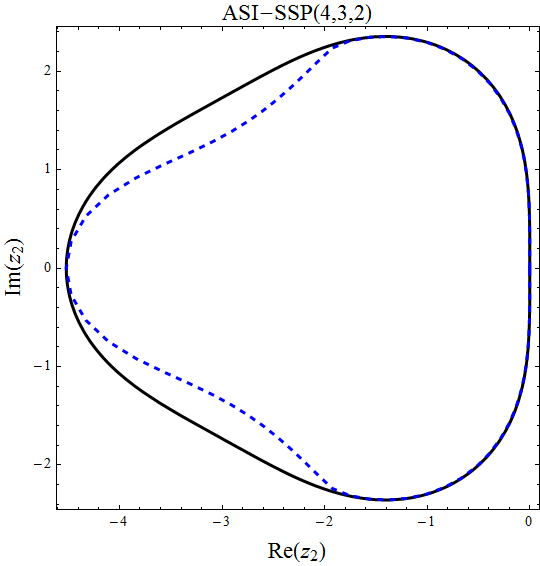}
  \includegraphics[width=0.4\textwidth]{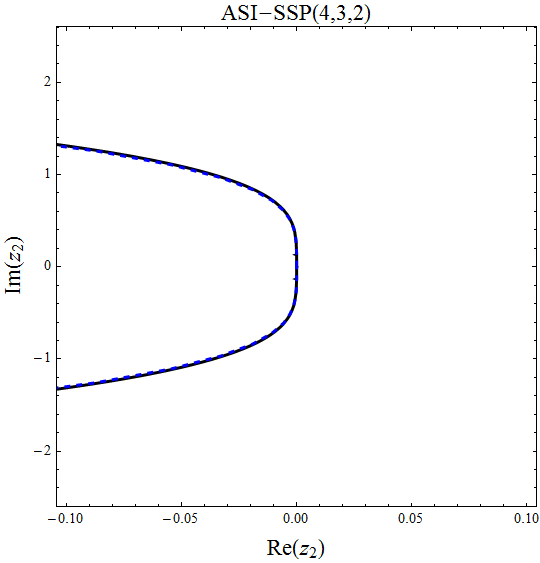}\\
  \includegraphics[width=0.4\textwidth]{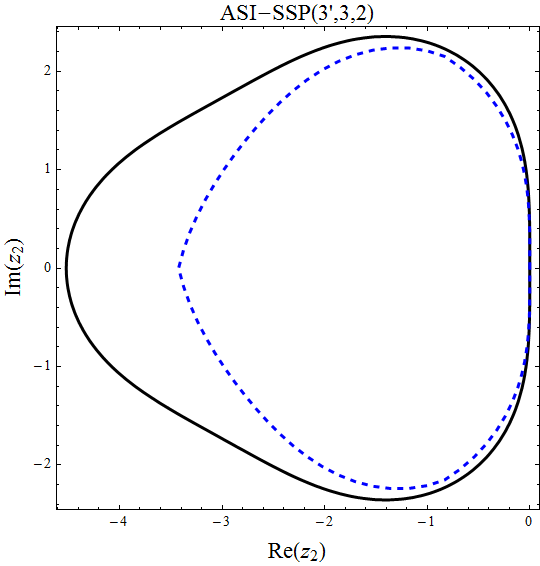}
  \includegraphics[width=0.4\textwidth]{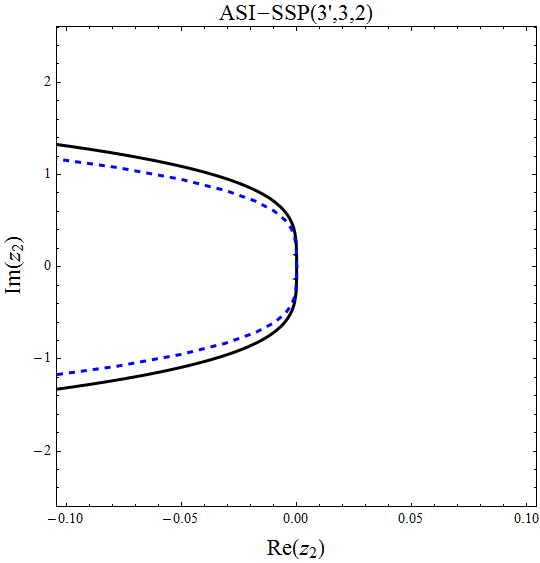}\\
  \includegraphics[width=0.4\textwidth]{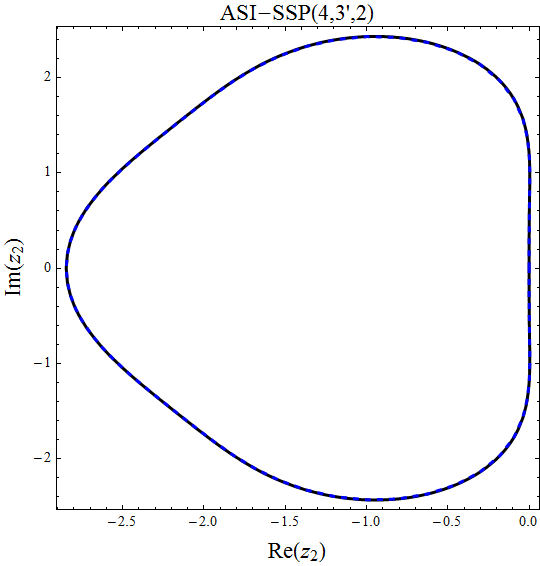}
  \includegraphics[width=0.4\textwidth]{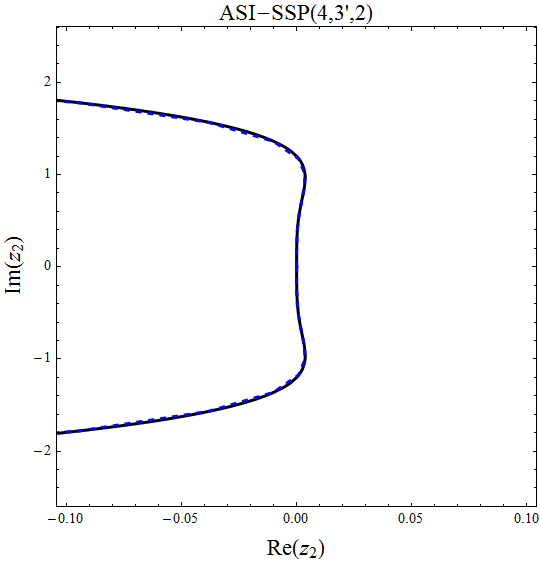}\\
\end{figure}
\begin{figure}[htbp]
\centering
  \includegraphics[width=0.4\textwidth]{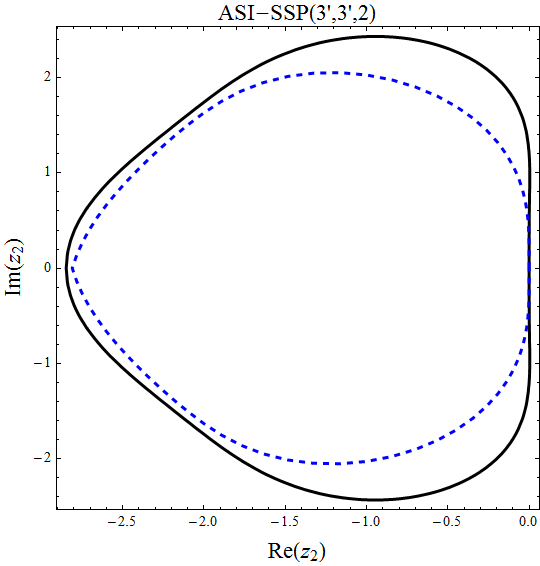}
  \includegraphics[width=0.4\textwidth]{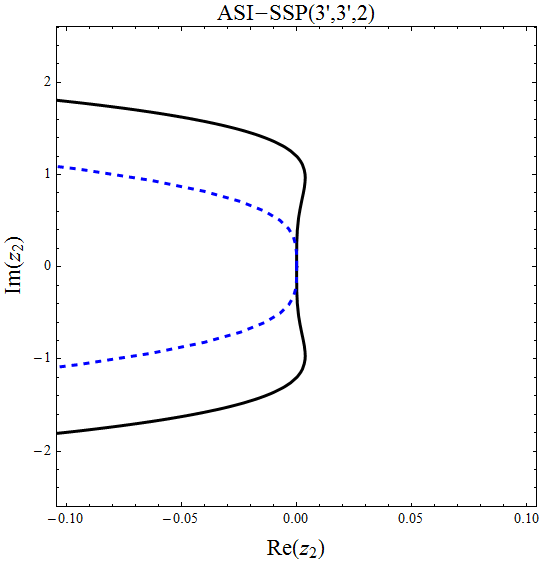}\\
  \includegraphics[width=0.4\textwidth]{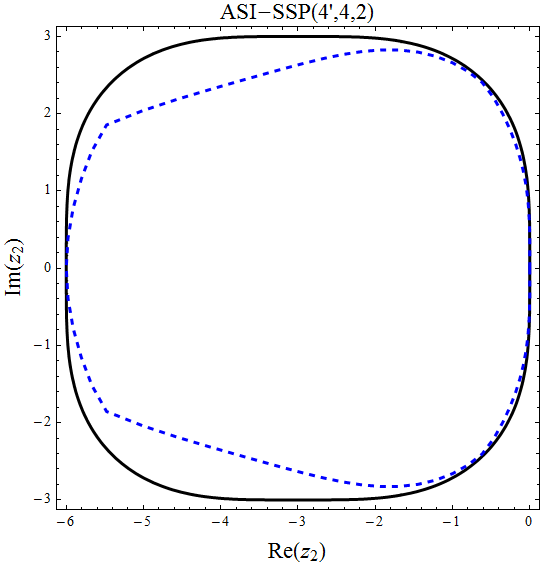}
  \includegraphics[width=0.4\textwidth]{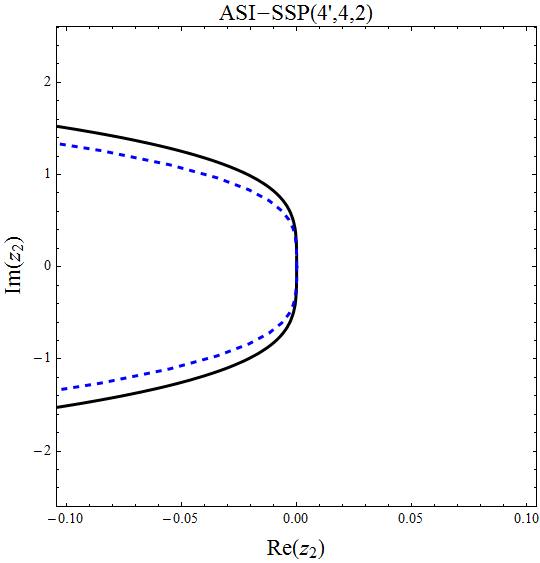}\\
  \includegraphics[width=0.4\textwidth]{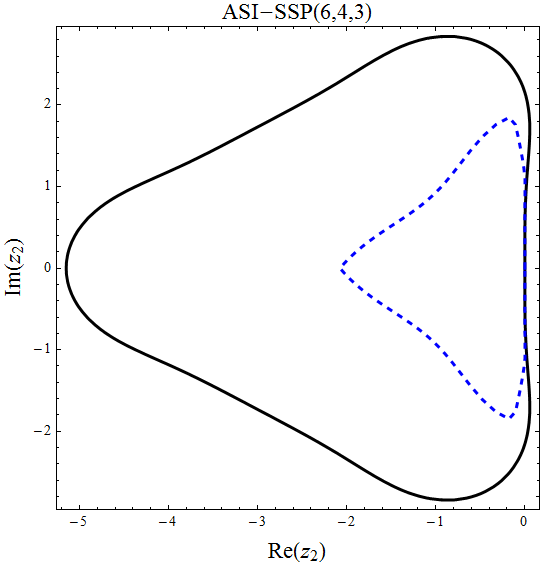}
  \includegraphics[width=0.4\textwidth]{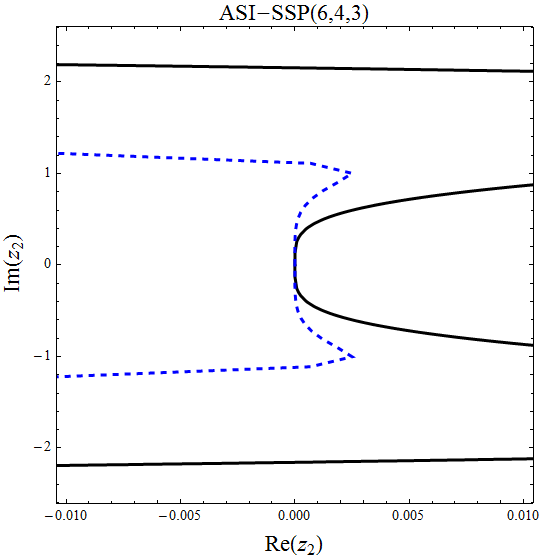}\\
\end{figure}
\begin{figure}[htbp]
\centering
  \includegraphics[width=0.4\textwidth]{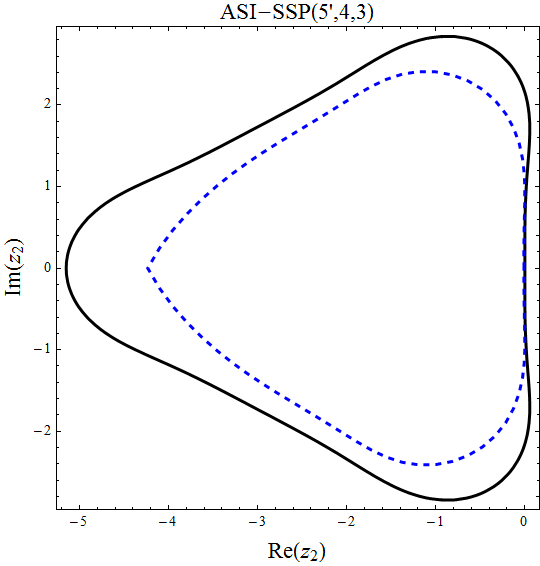}
  \includegraphics[width=0.4\textwidth]{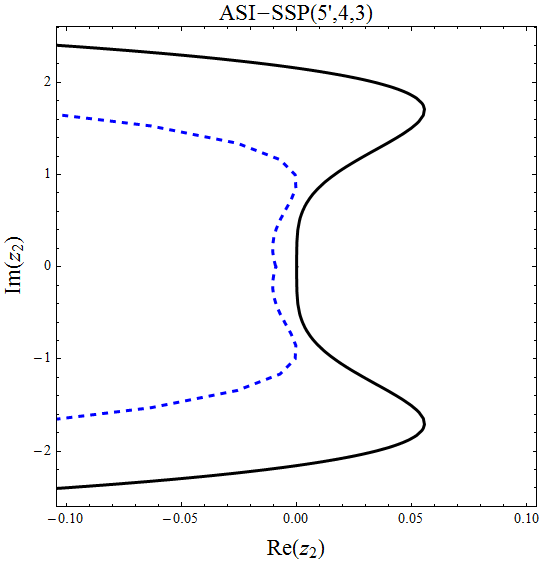}\\
  \includegraphics[width=0.4\textwidth]{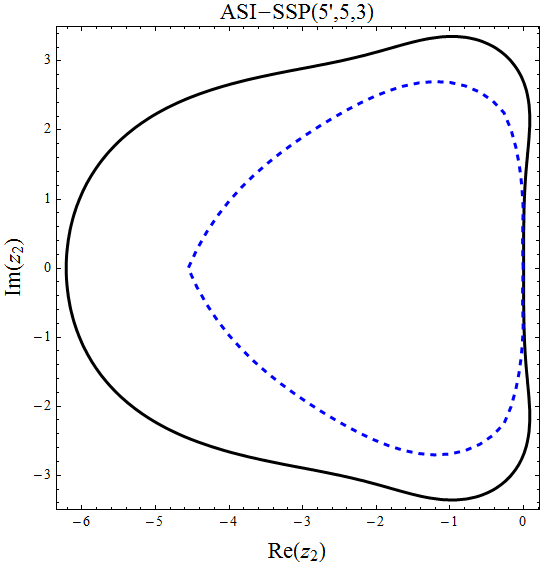}
  \includegraphics[width=0.4\textwidth]{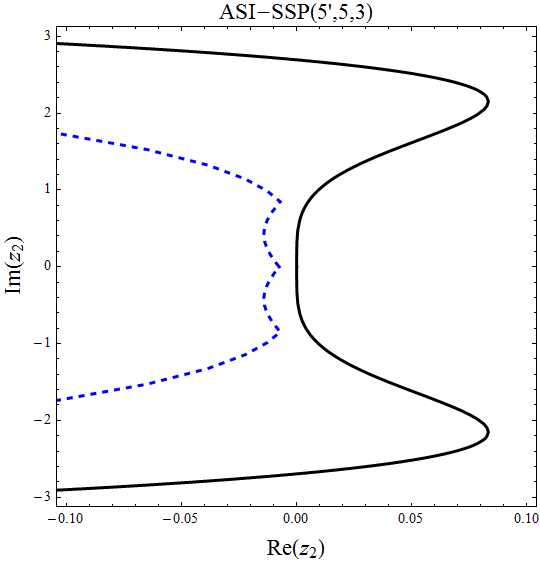}\\
  \includegraphics[width=0.4\textwidth]{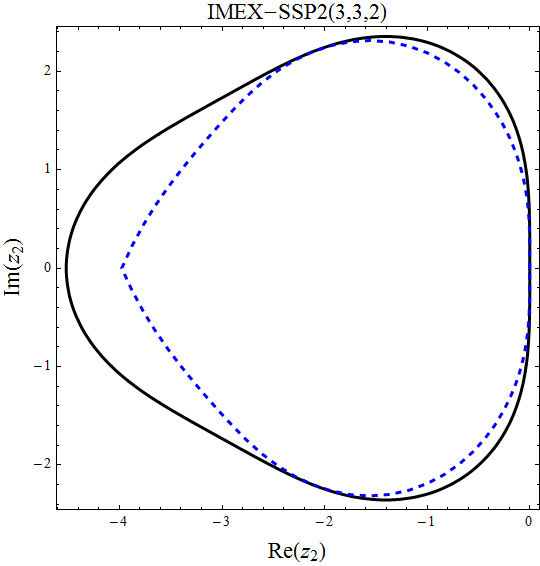}
  \includegraphics[width=0.4\textwidth]{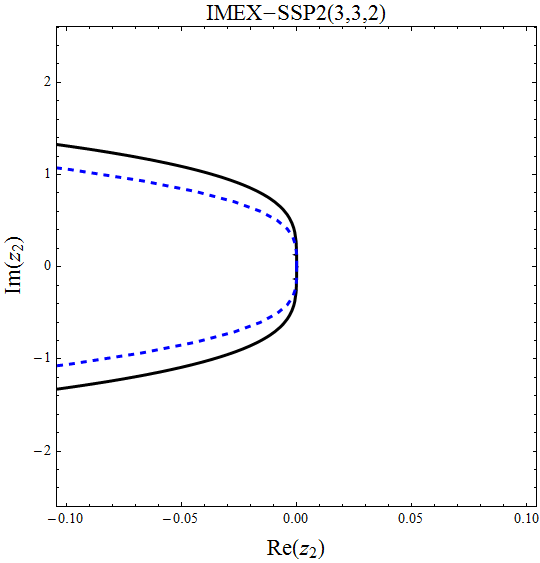}\\
\end{figure}

\section{Tests}

We first consider Pareschi and Russo's problem \cite{Pareschi:2001} to study the convergence behaviors in the temporal domain:
\begin{equation}
\dot{x}(t) = -y(t), \dot{y}(t) = (\sin(x(t))-y(t))/\varepsilon + x(t).
\end{equation}
The eigenvalues of the explicit part are $\pm i$. To partition for an IMEX scheme, the terms divided by $\varepsilon$ are integrated with the implicit method, whereas the other terms are integrated explicitly. The initial conditions are considered in two forms: equilibrium initial conditions accomplished with 
$x(0)=\pi/2, y(0)=1$, and non-equilibrium, or perturbed, conditions specified by replacing the condition on $y$ with $y(0)=1/2$.

The $L_2$-norm error over the numerical integration interval $t\in[0,5]$ is used to assess the convergence behaviors. The $L_2$-norm is a type of average defined \cite{Liu:2006} as
\begin{equation}
E(\Delta t) = \sqrt{\Delta t \sum_i | x^i - x(t_i) |^2 },
\end{equation}
where $x^i$ is computed by the IMEX scheme and $x(t_i)$ is expected to be an exact solution, but is substituted by the computed solution on the finest time grids. The error is a function of $\varepsilon$ and $\Delta t$, that is,  $E=E(\varepsilon, \Delta t)$. We examine it in two-dimensional parameter space $(\varepsilon, \Delta t)$.

\begin{figure}[htbp]
\centering
  \includegraphics[width=0.45\textwidth]{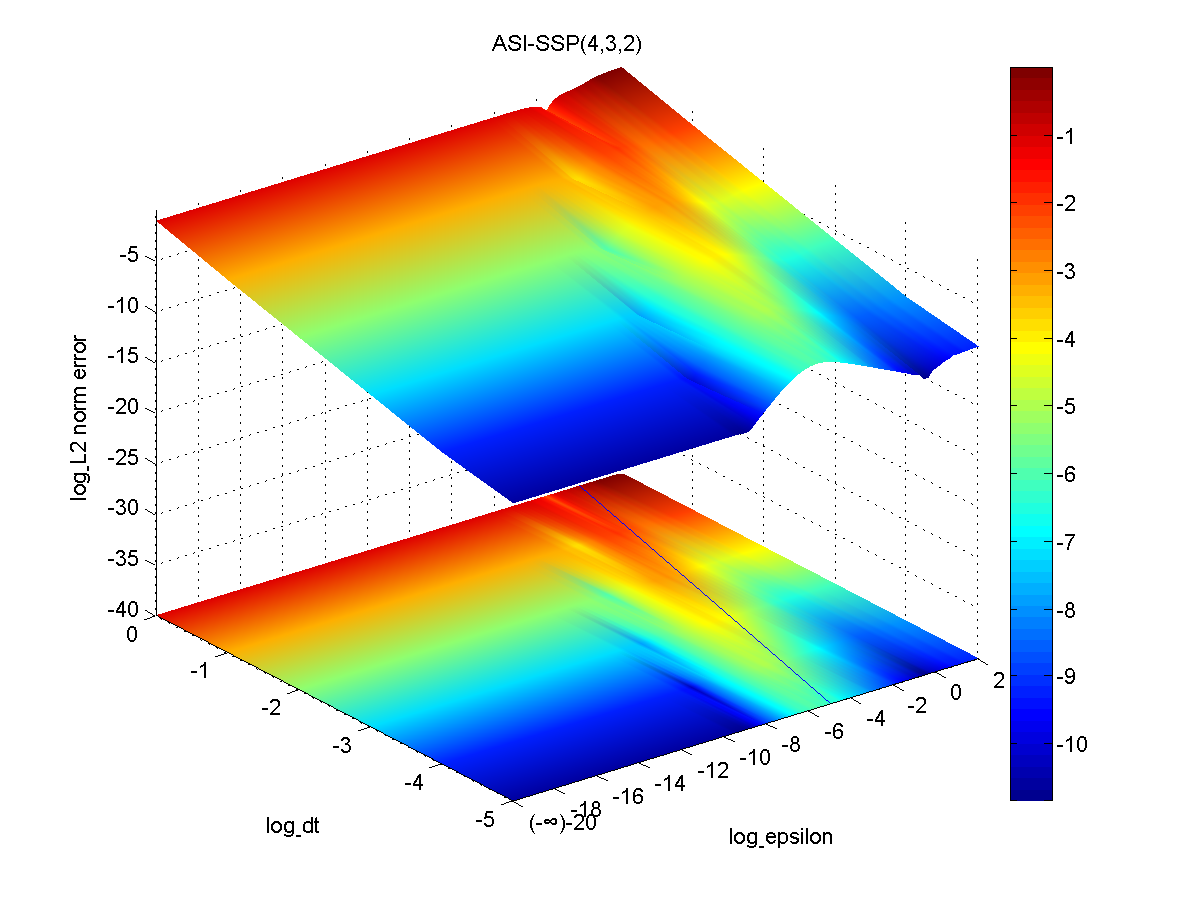}
  \includegraphics[width=0.45\textwidth]{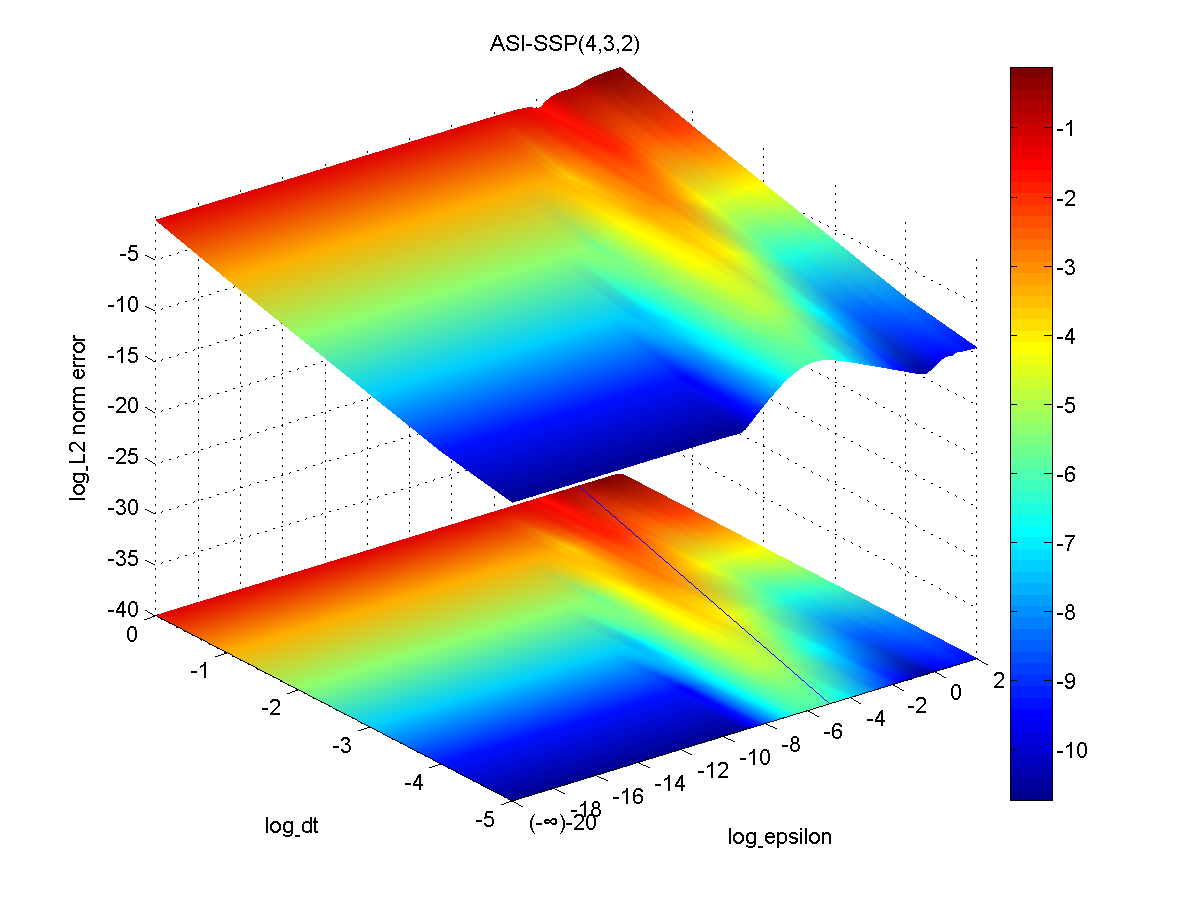}\\
  \includegraphics[width=0.45\textwidth]{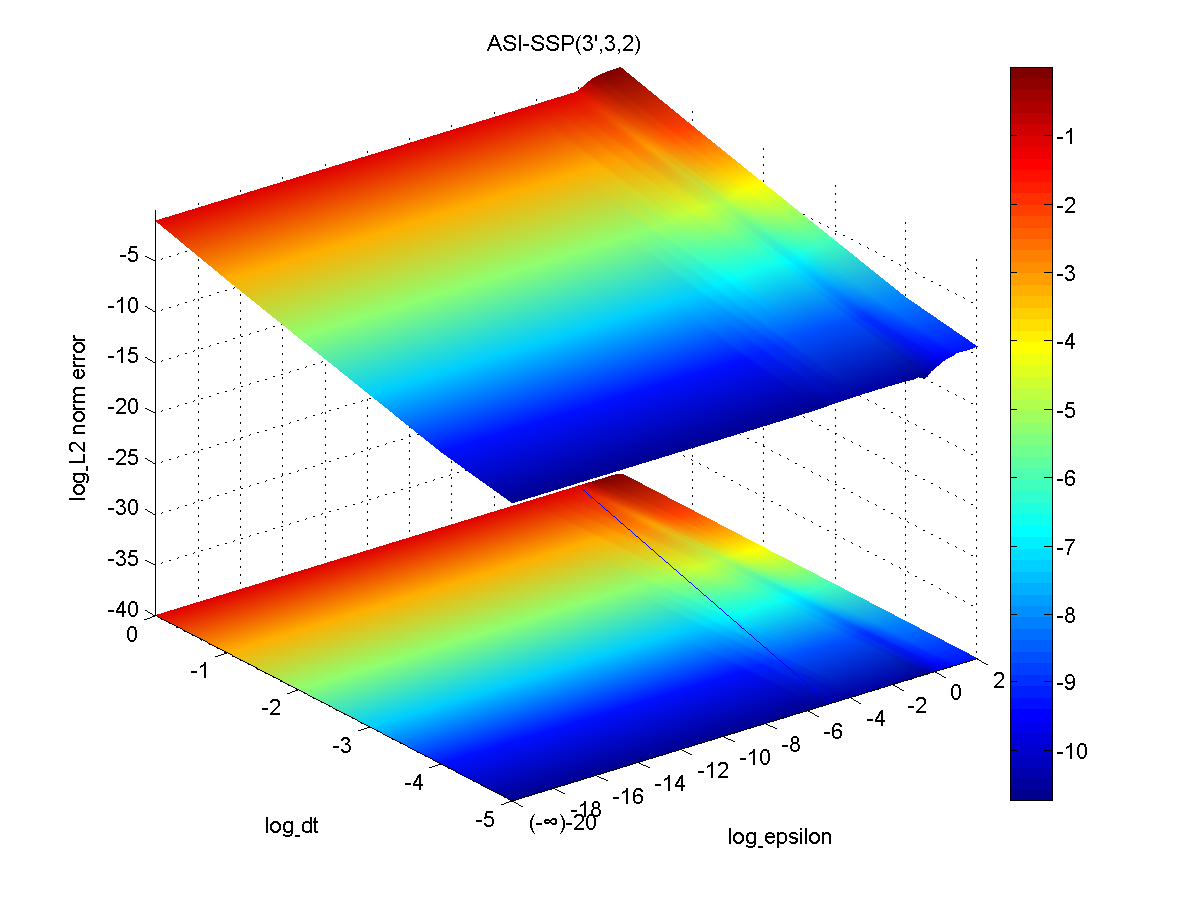}
  \includegraphics[width=0.45\textwidth]{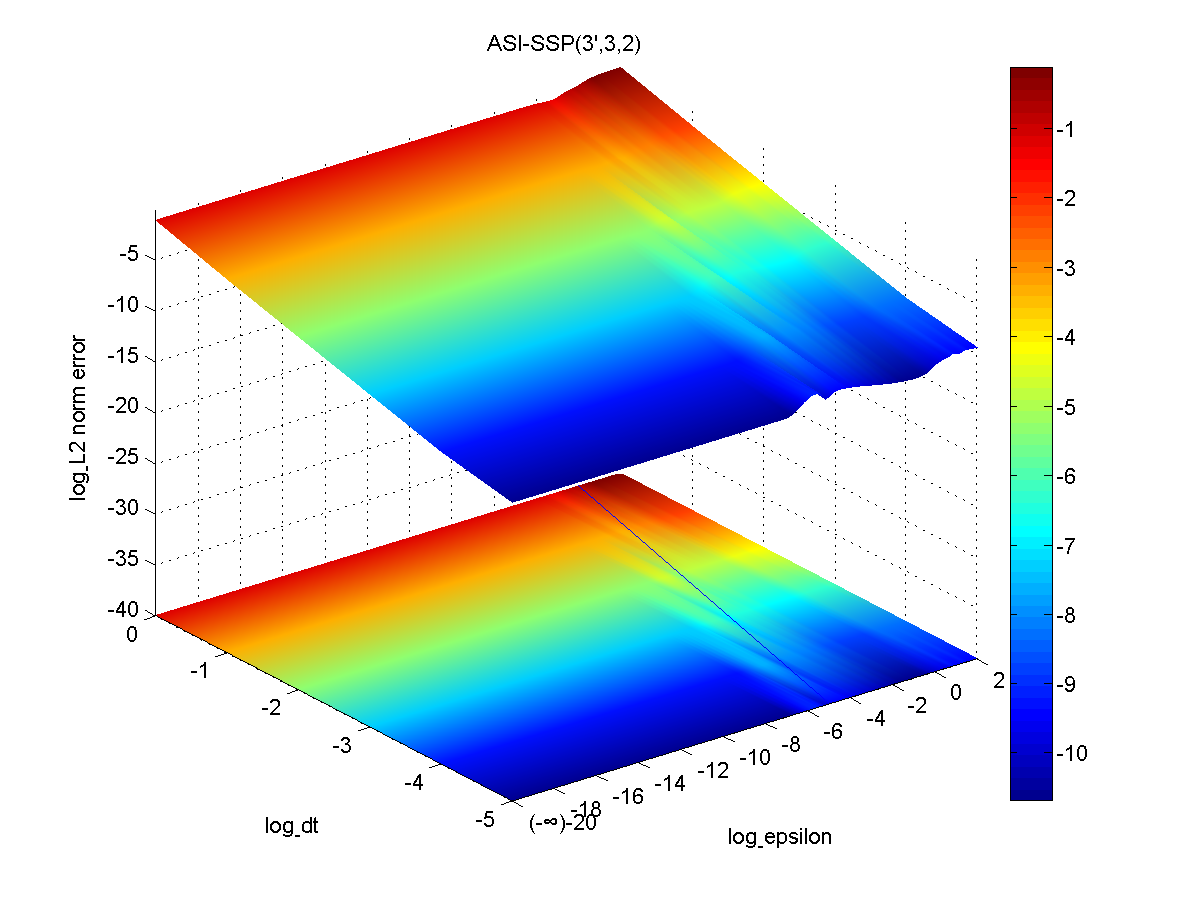}\\
\caption{Representative multi-dimensional displays for the $L_2$-norm error tested with two selected schemes for Pareschi and Russo's problem with equilibrium initial conditions. The left column is for $x$ and the right column is for $y$.}\label{fig.L2NormError}
\end{figure}

\begin{figure}[htbp]
\centering
  \includegraphics[width=0.38\textwidth]{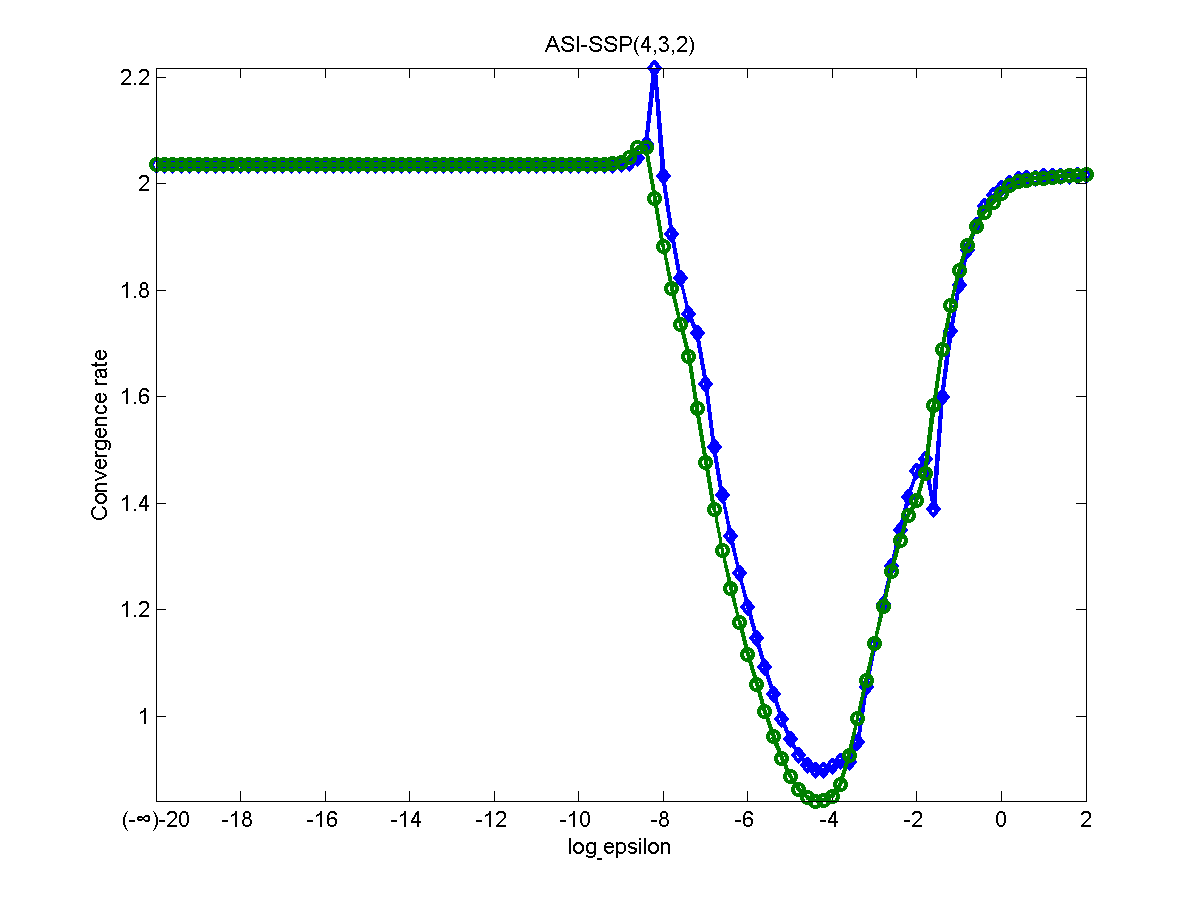}
  \includegraphics[width=0.38\textwidth]{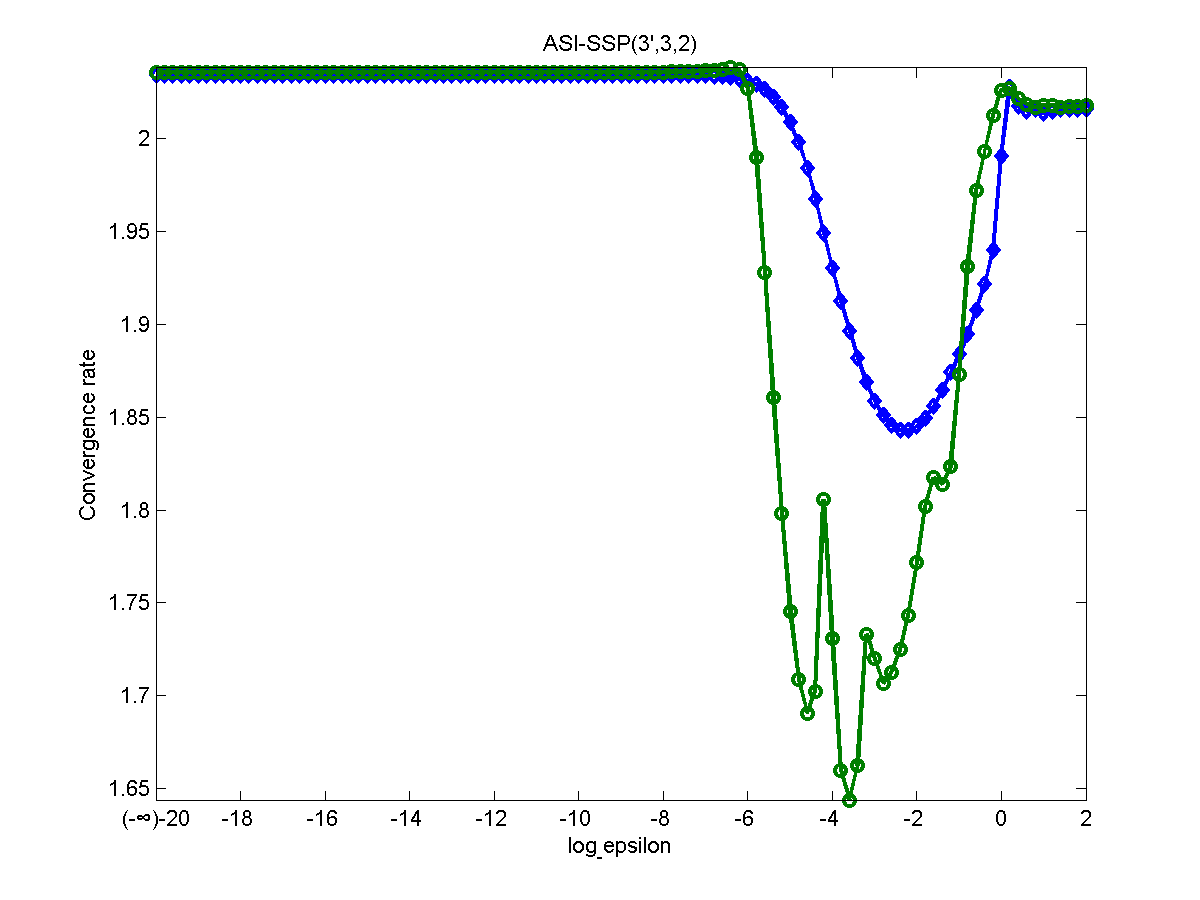}\\
  \includegraphics[width=0.38\textwidth]{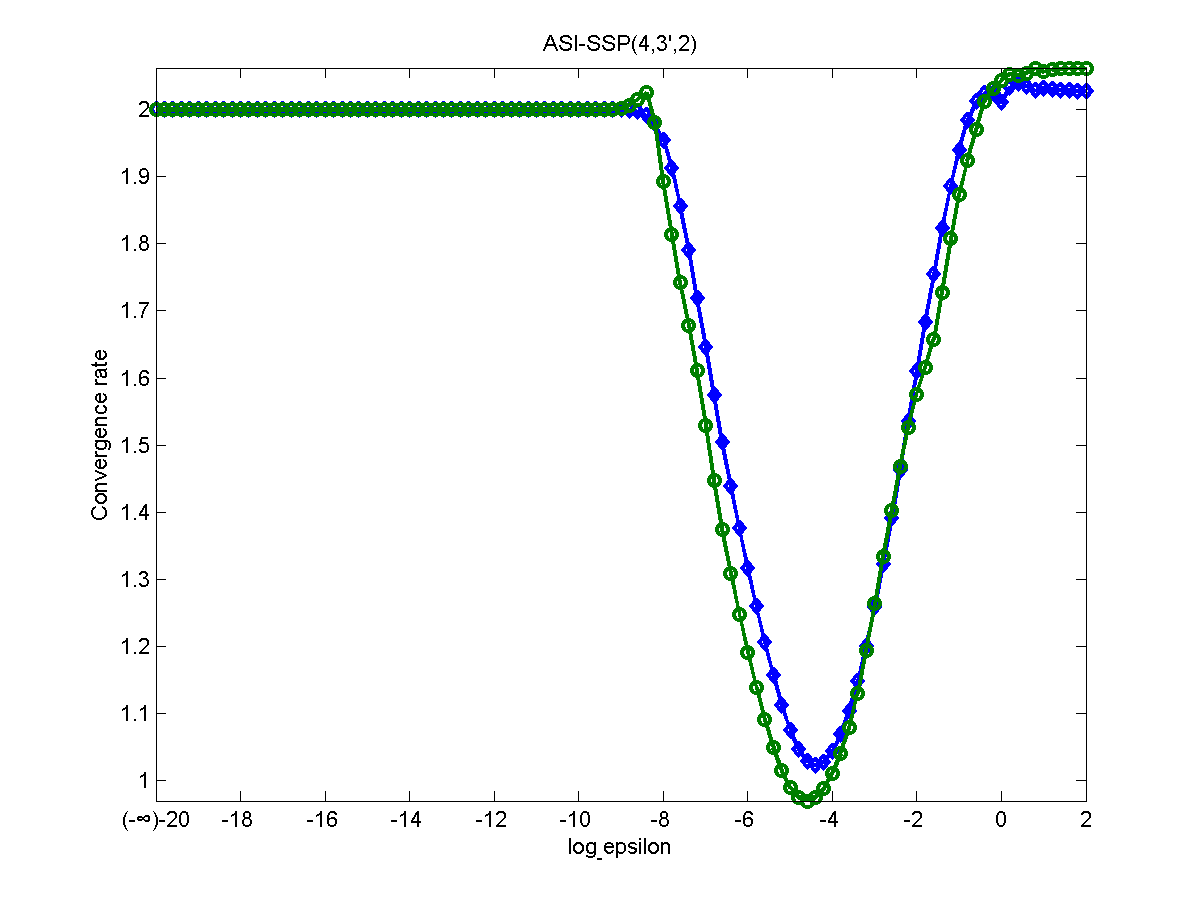}
  \includegraphics[width=0.38\textwidth]{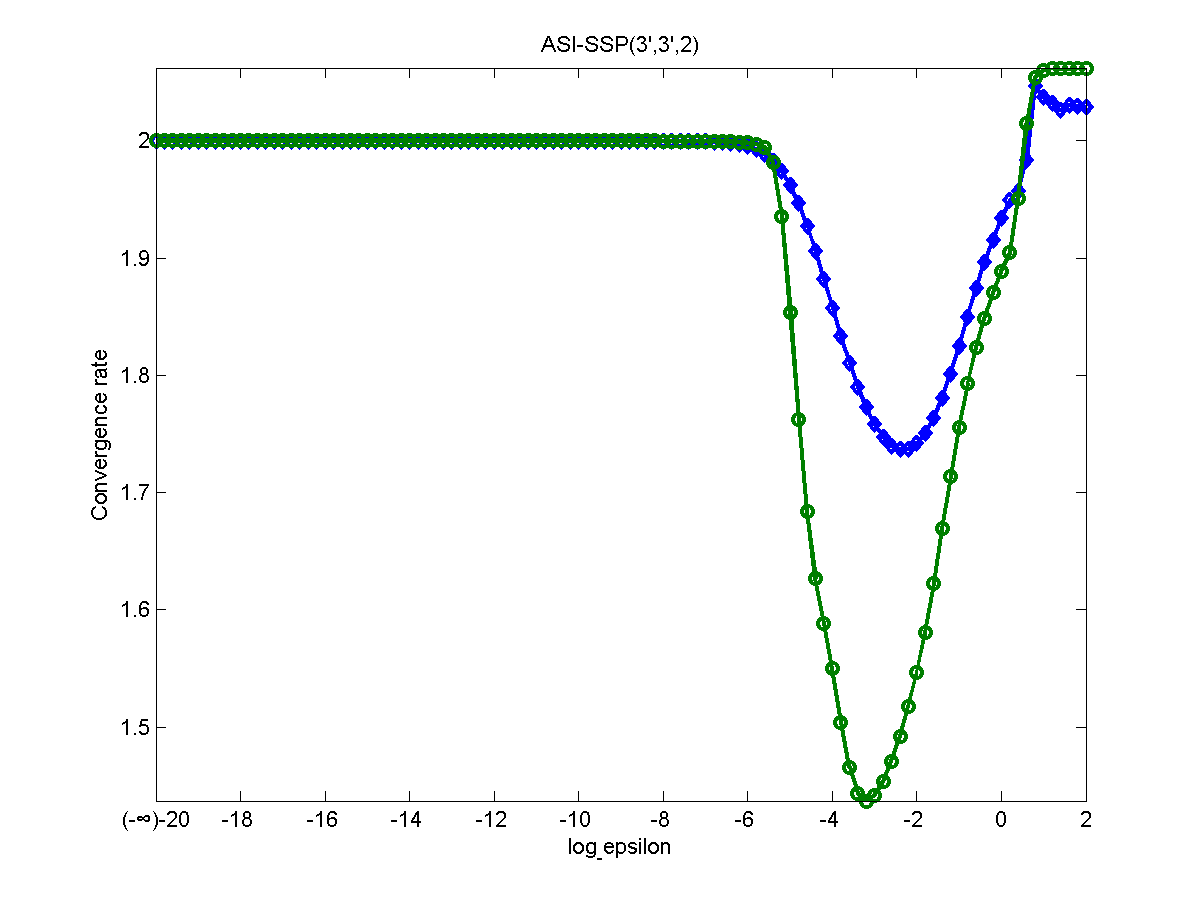}\\
  \includegraphics[width=0.38\textwidth]{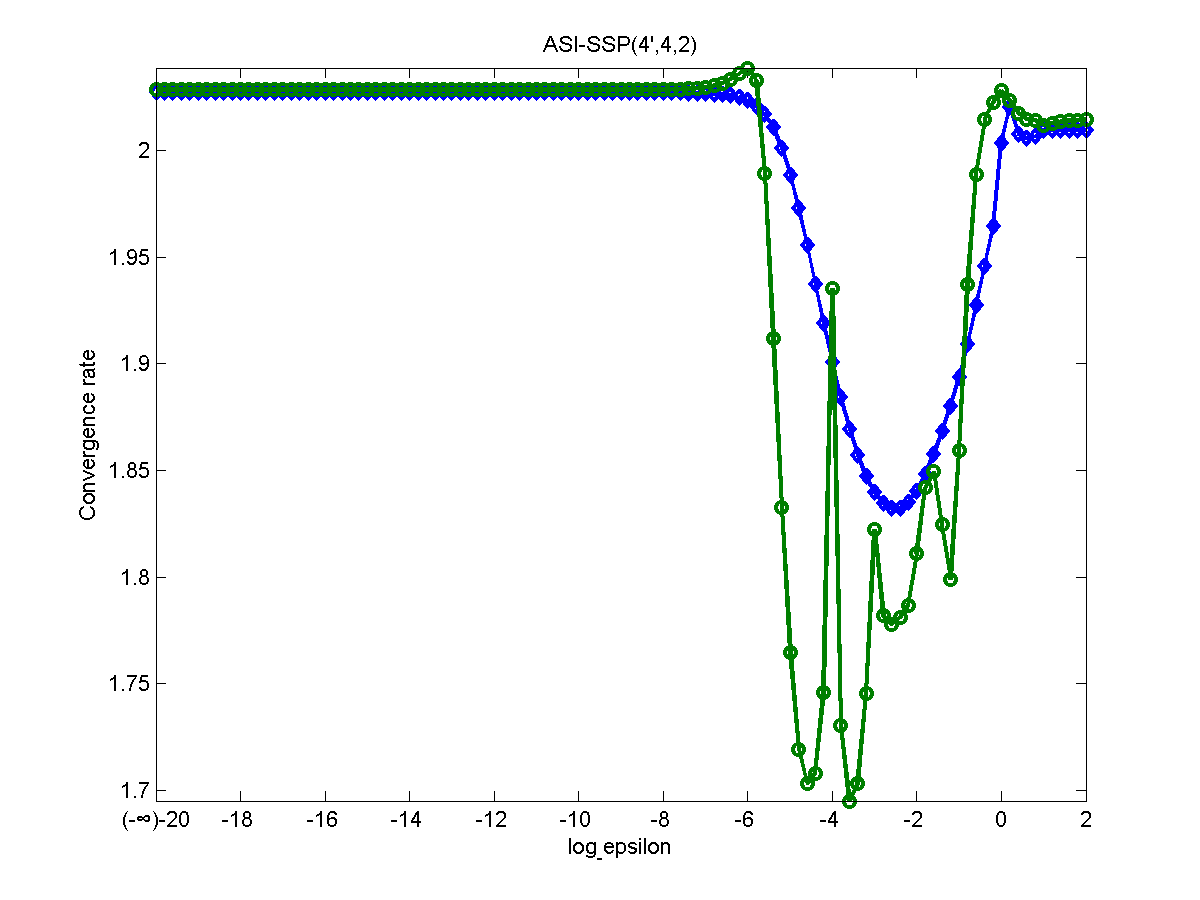}
  \includegraphics[width=0.38\textwidth]{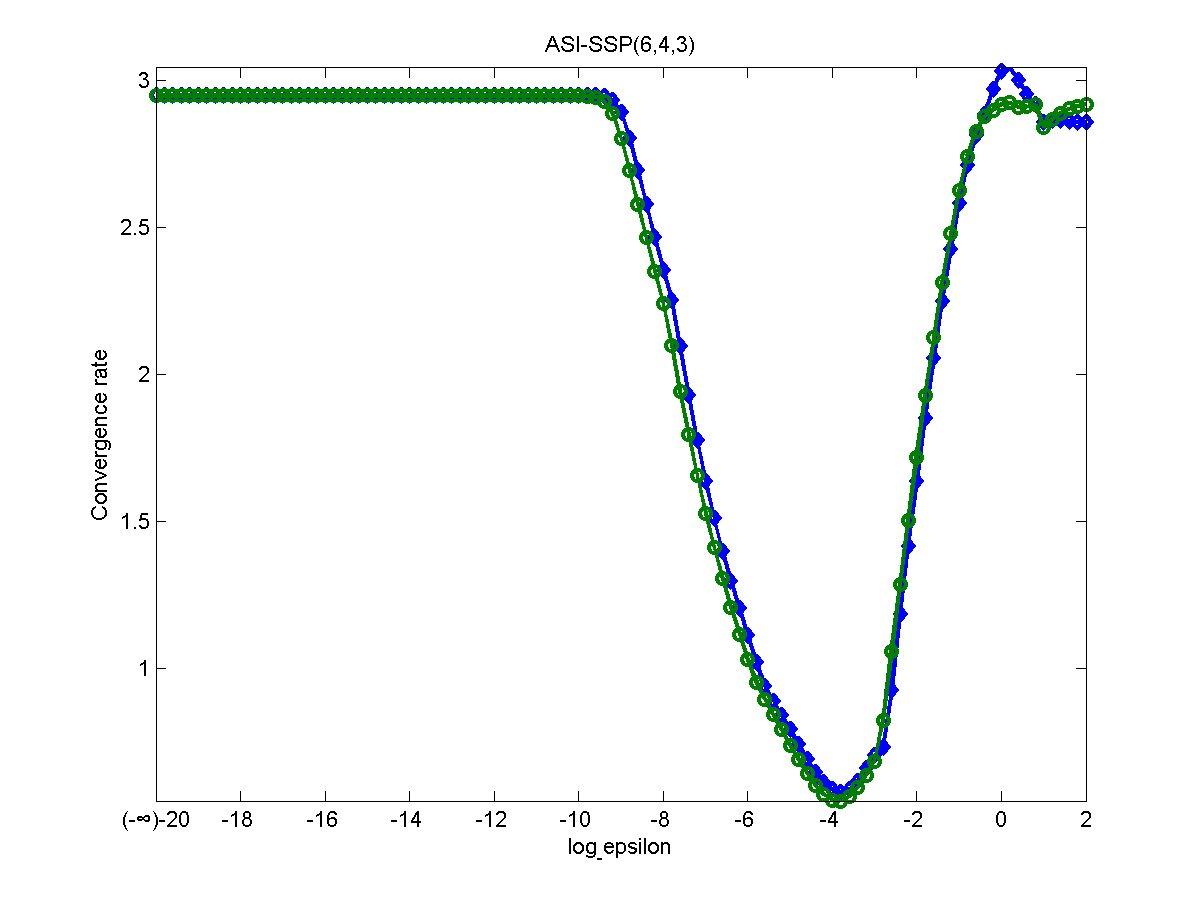}\\
  \includegraphics[width=0.38\textwidth]{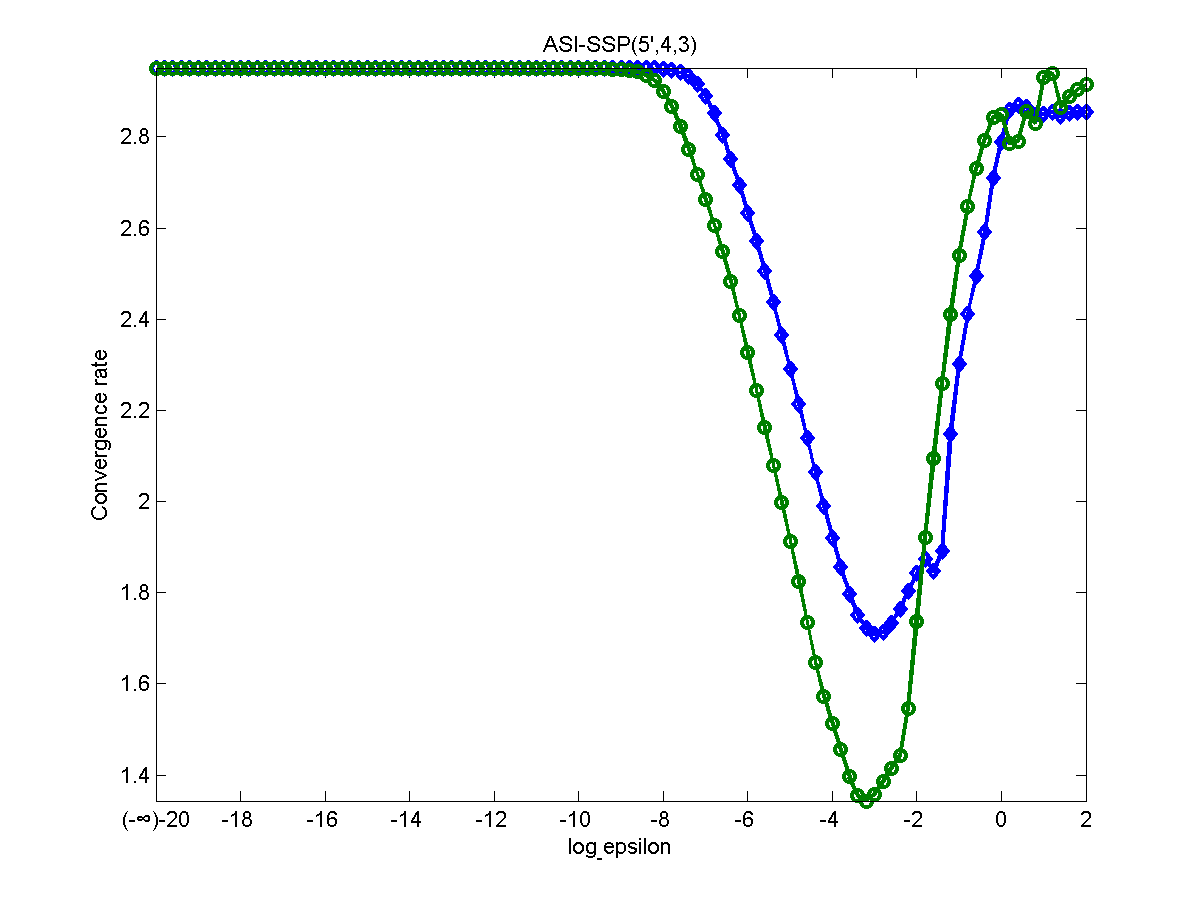}
  \includegraphics[width=0.38\textwidth]{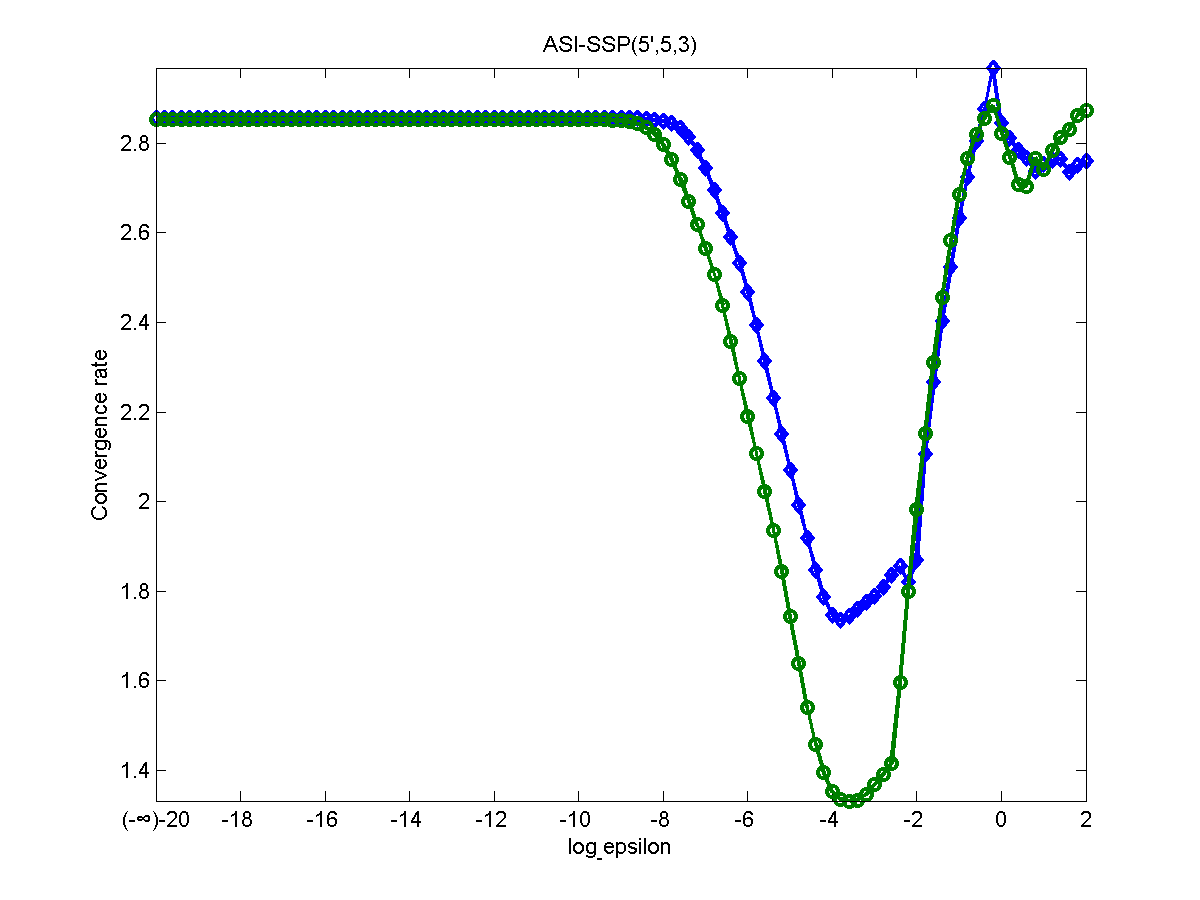}\\
  \includegraphics[width=0.38\textwidth]{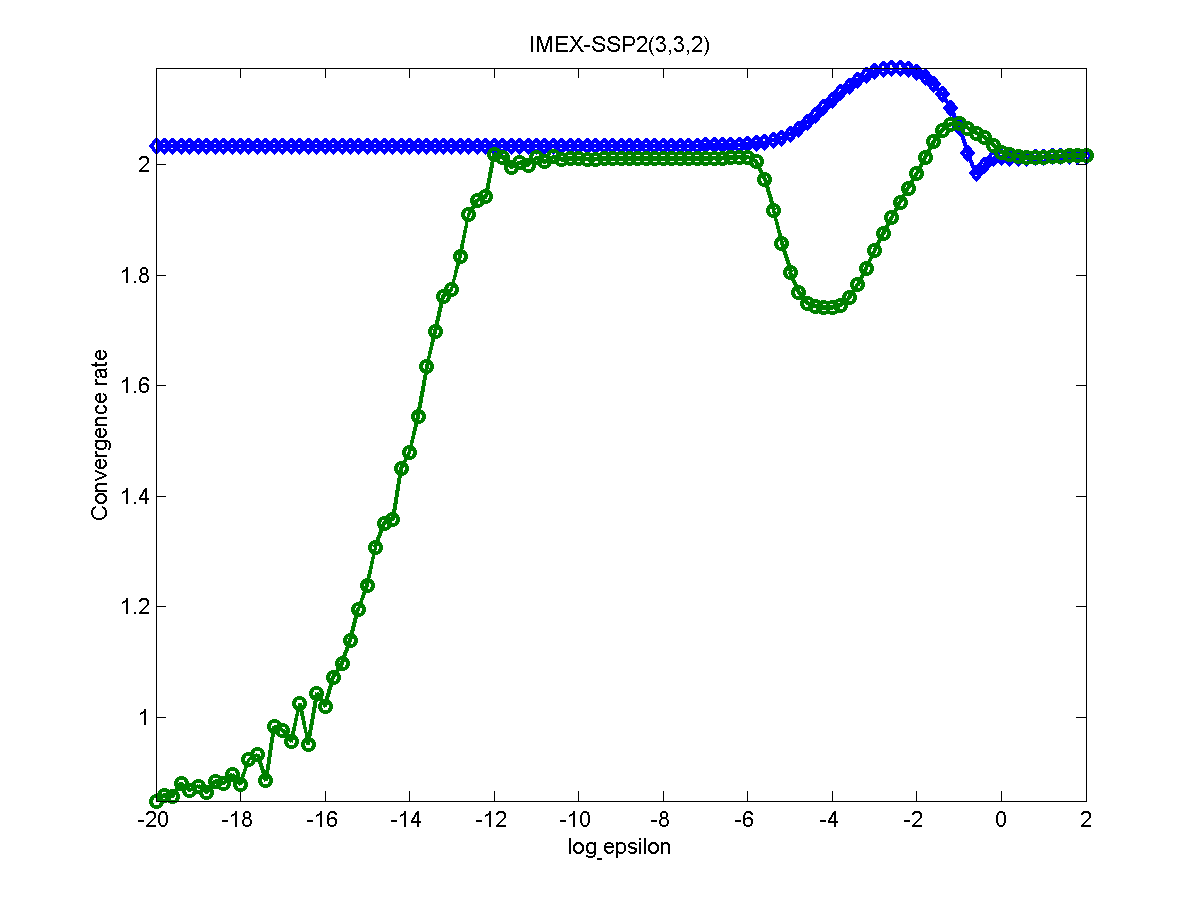}
  \includegraphics[width=0.38\textwidth]{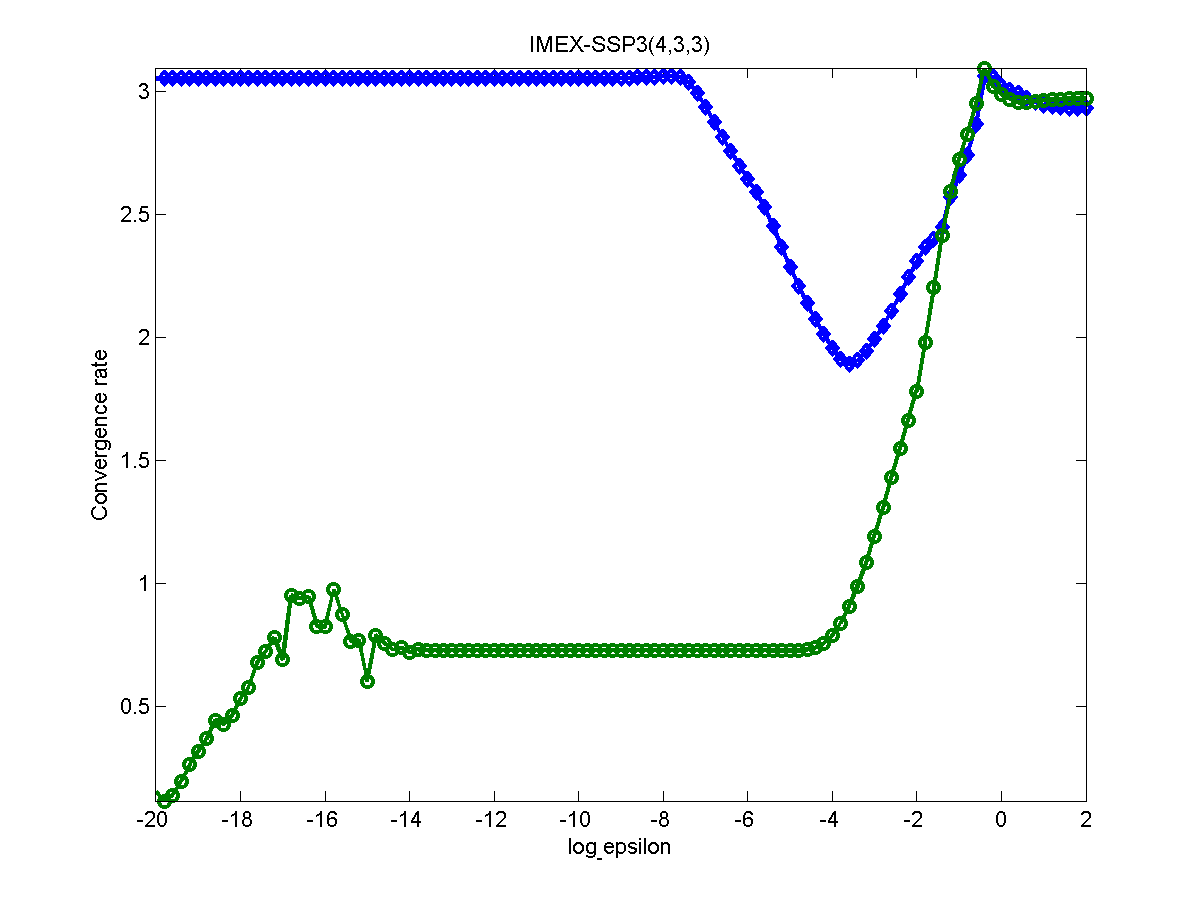}\\
\caption{Convergence behaviors of eight schemes designed in this paper and two from Pareschi and Russo \cite{Pareschi:2005} for Pareschi and Russo's problem with equilibrium initial conditions. The blue lines with diamonds are for $x$ and the green lines with circles are for $y$.}\label{fig.ConvergParWtEquil}
\end{figure}

\begin{figure}[htbp]
\centering
  \includegraphics[width=0.38\textwidth]{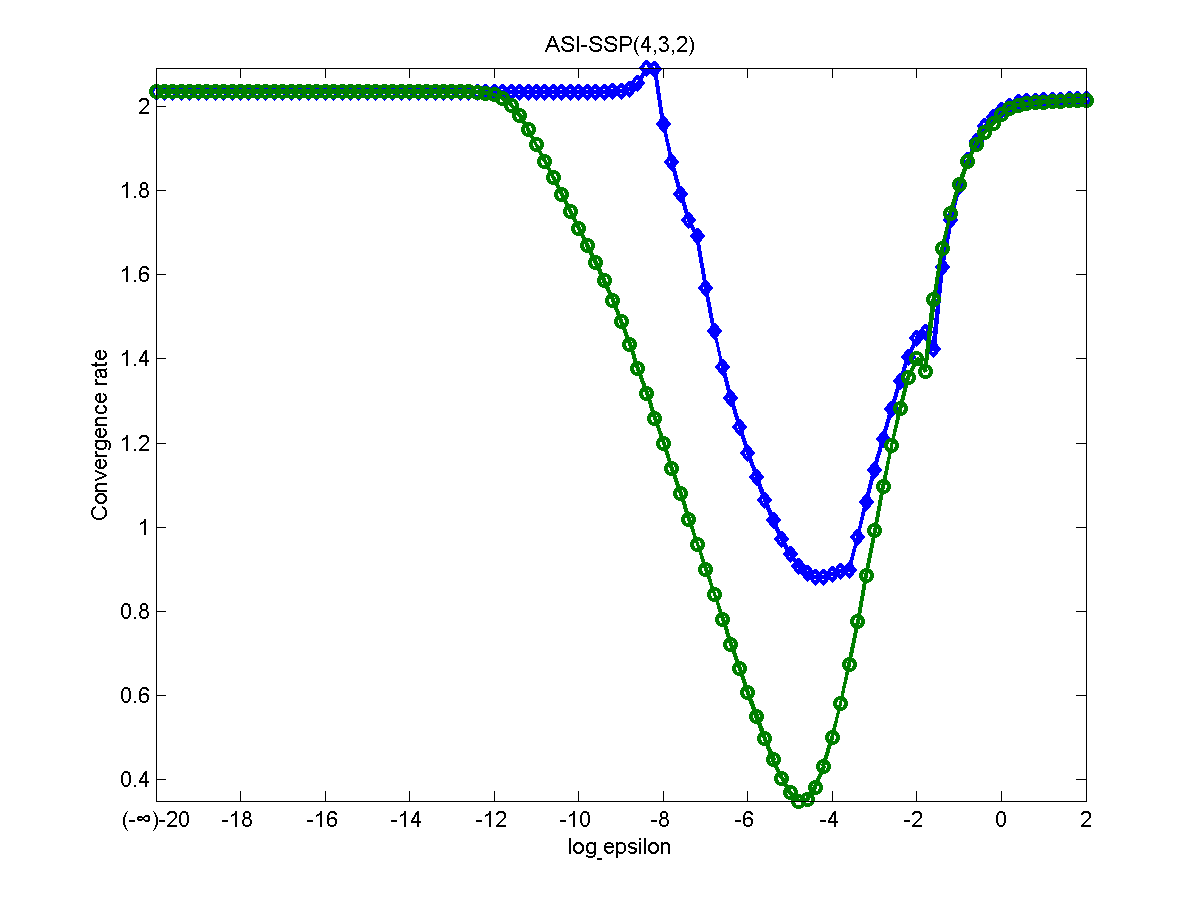}
  \includegraphics[width=0.38\textwidth]{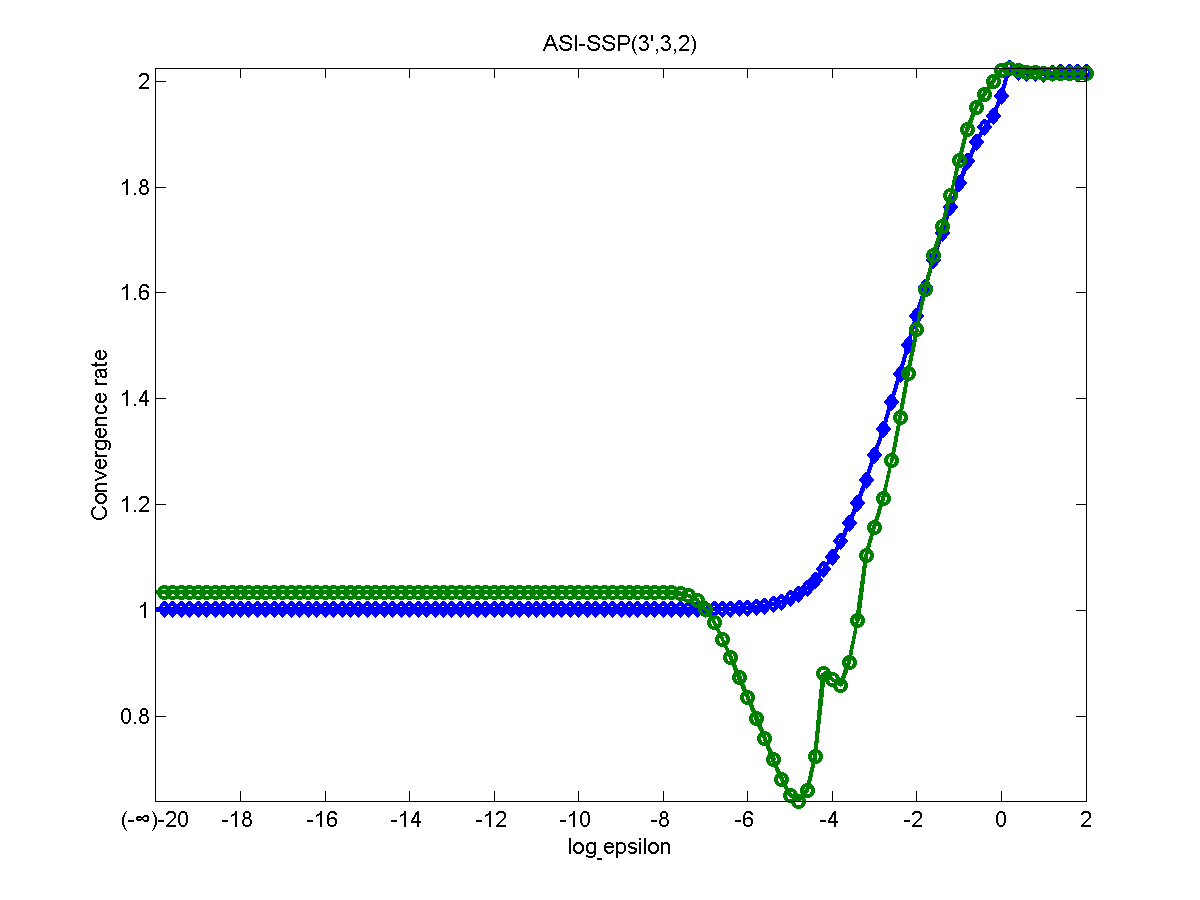}\\
  \includegraphics[width=0.38\textwidth]{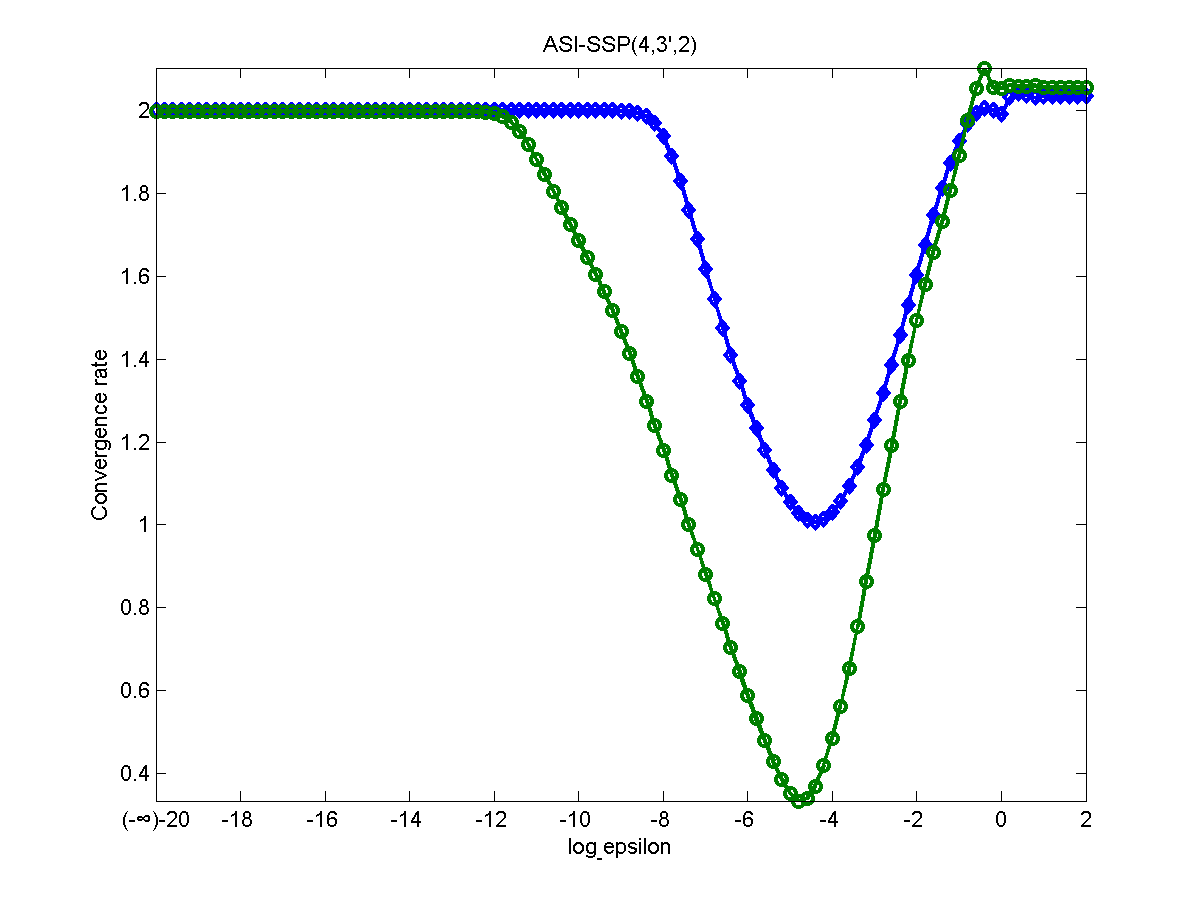}
  \includegraphics[width=0.38\textwidth]{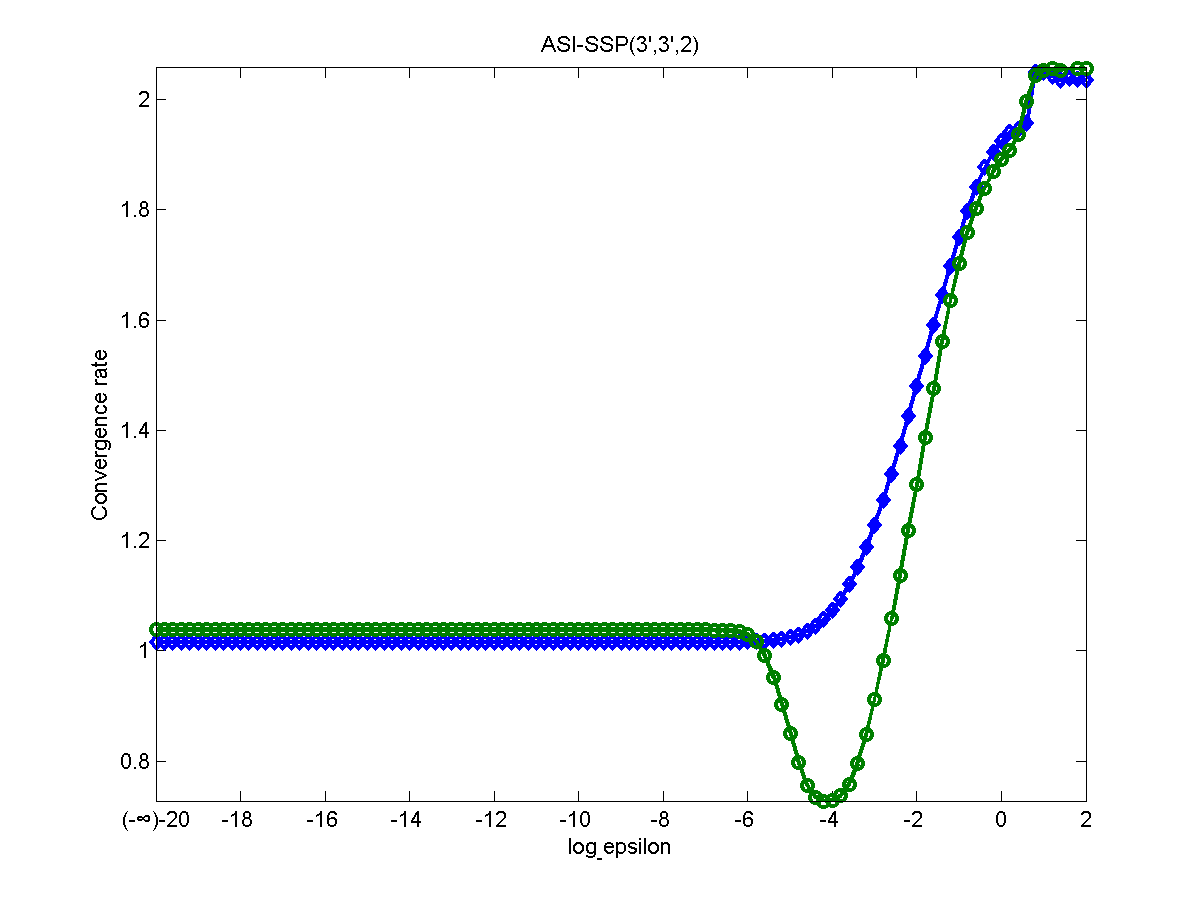}\\
  \includegraphics[width=0.38\textwidth]{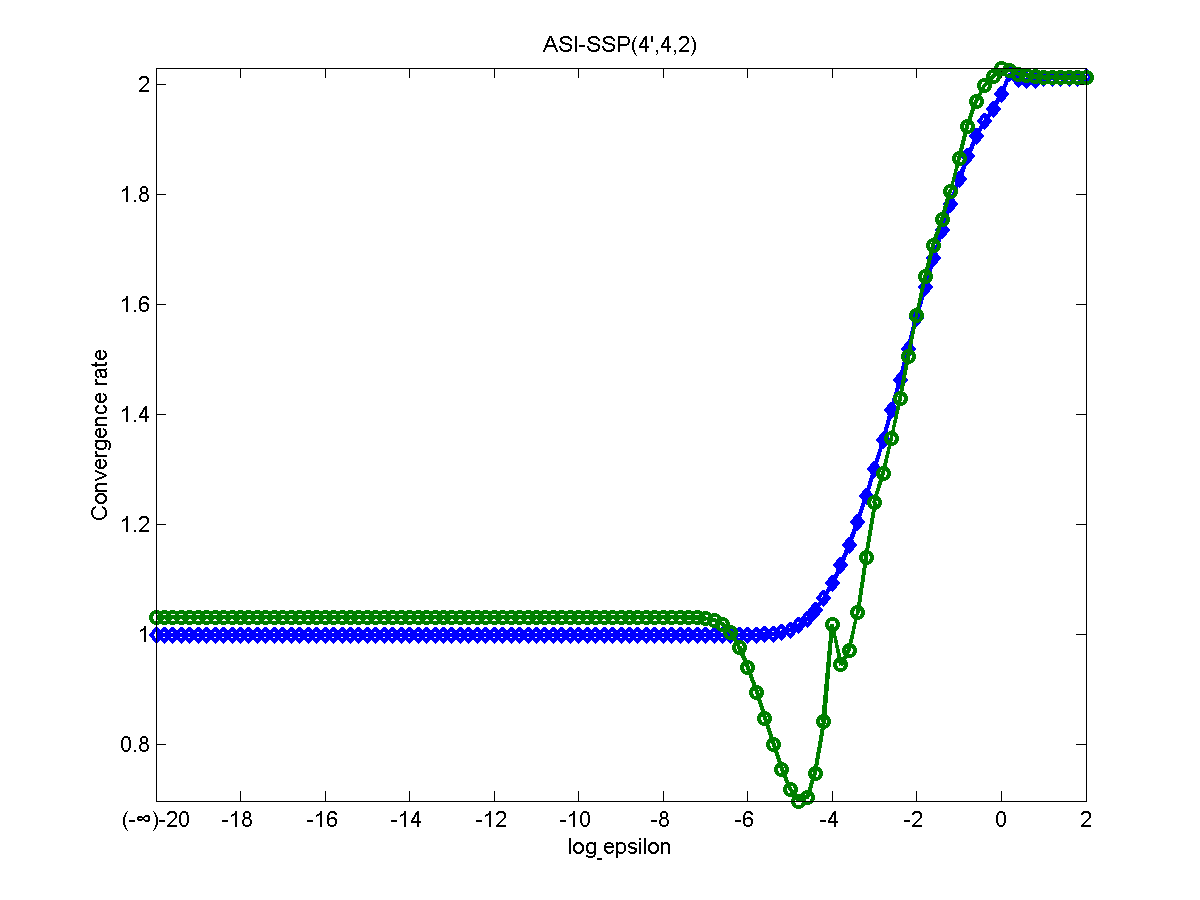}
  \includegraphics[width=0.38\textwidth]{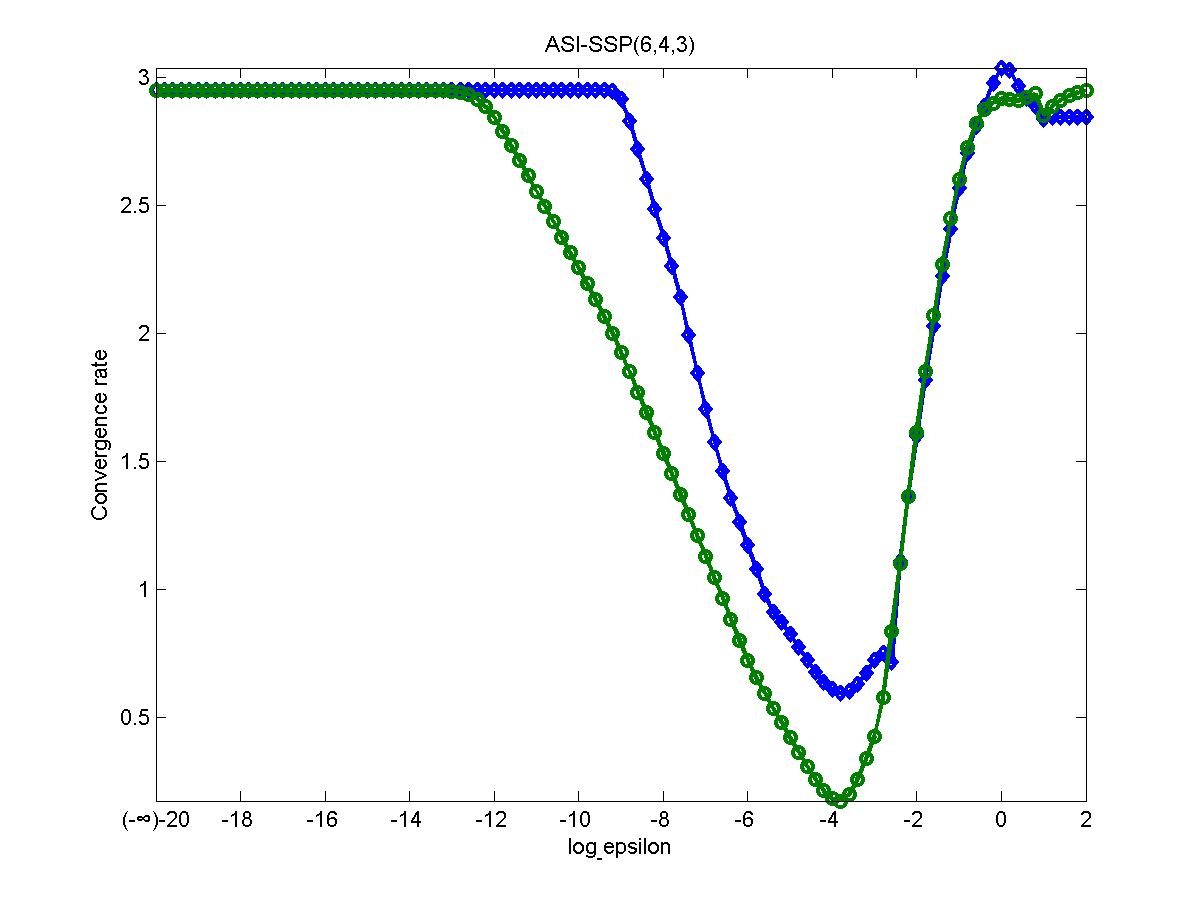}\\
  \includegraphics[width=0.38\textwidth]{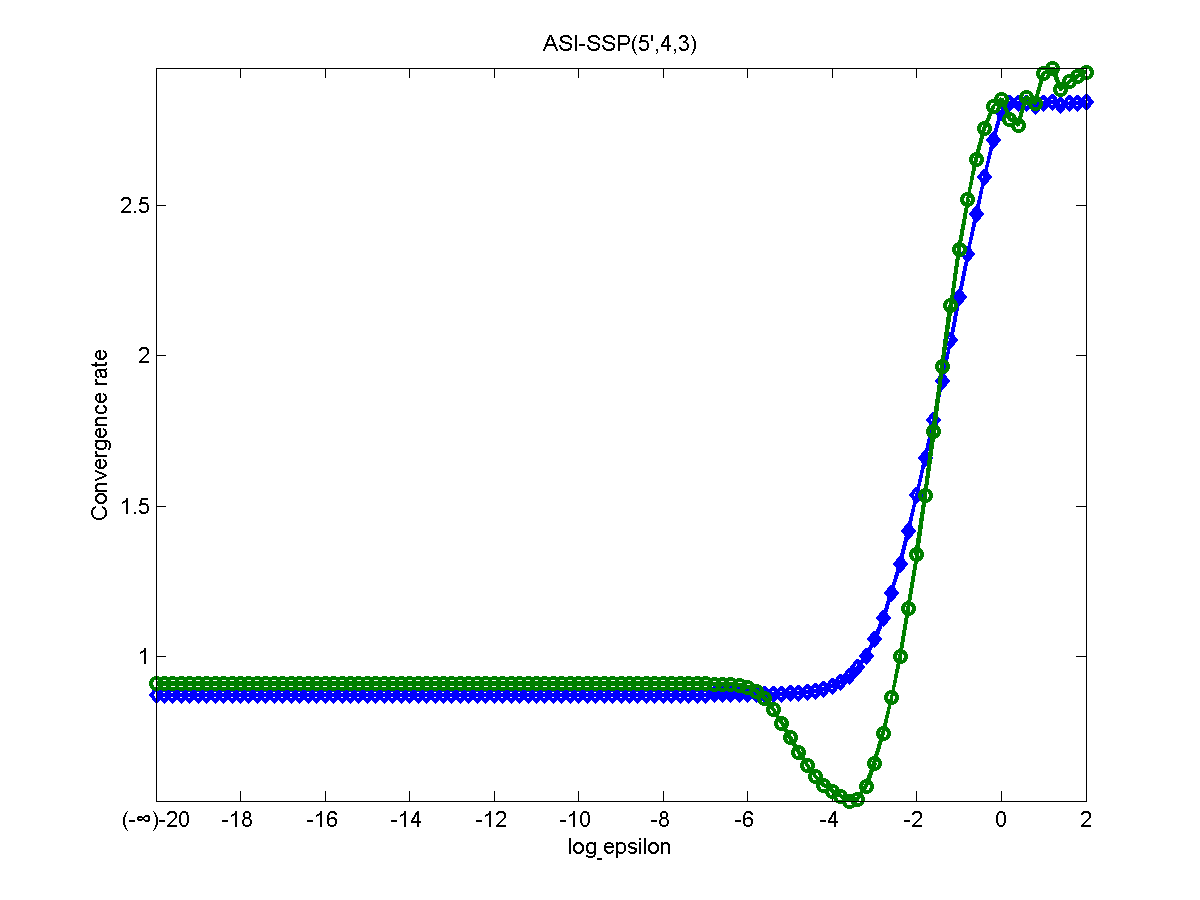}
  \includegraphics[width=0.38\textwidth]{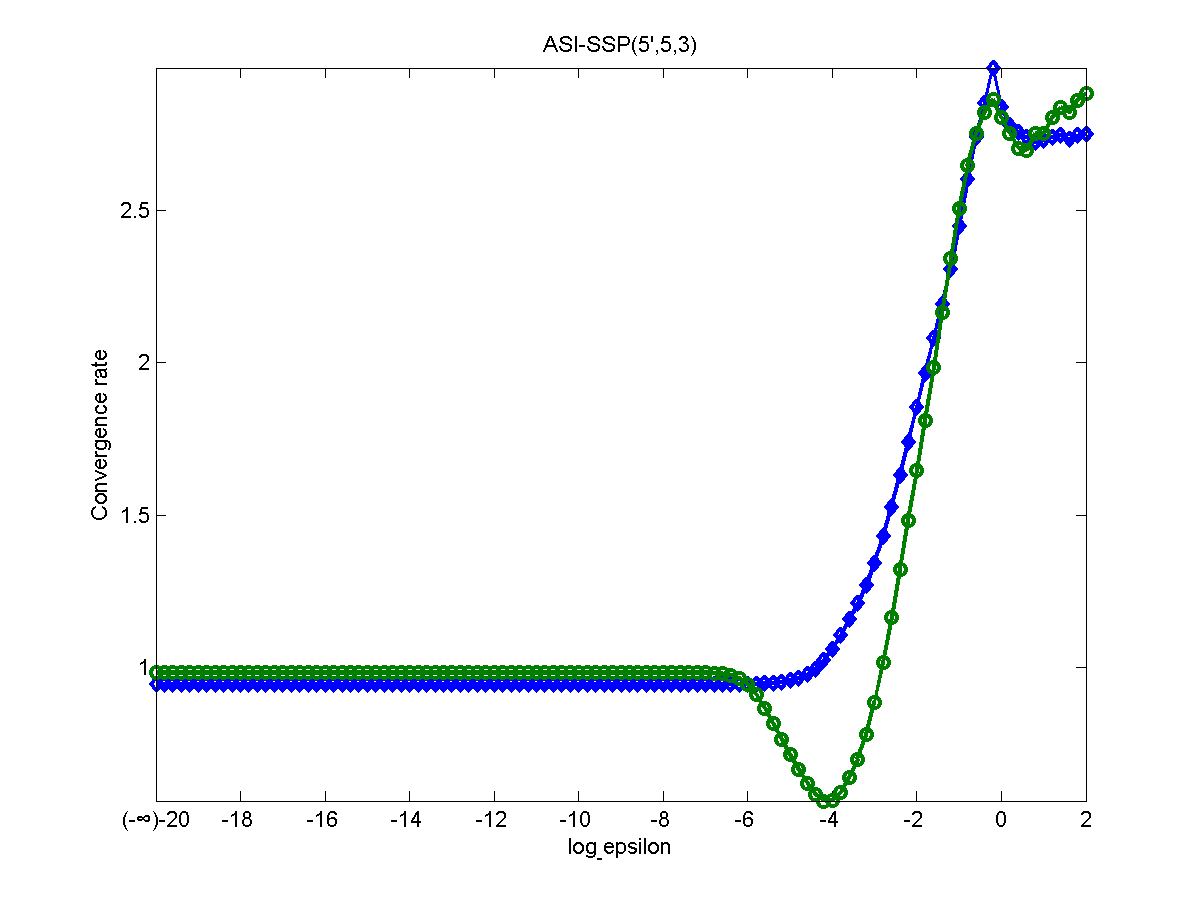}\\
  \includegraphics[width=0.38\textwidth]{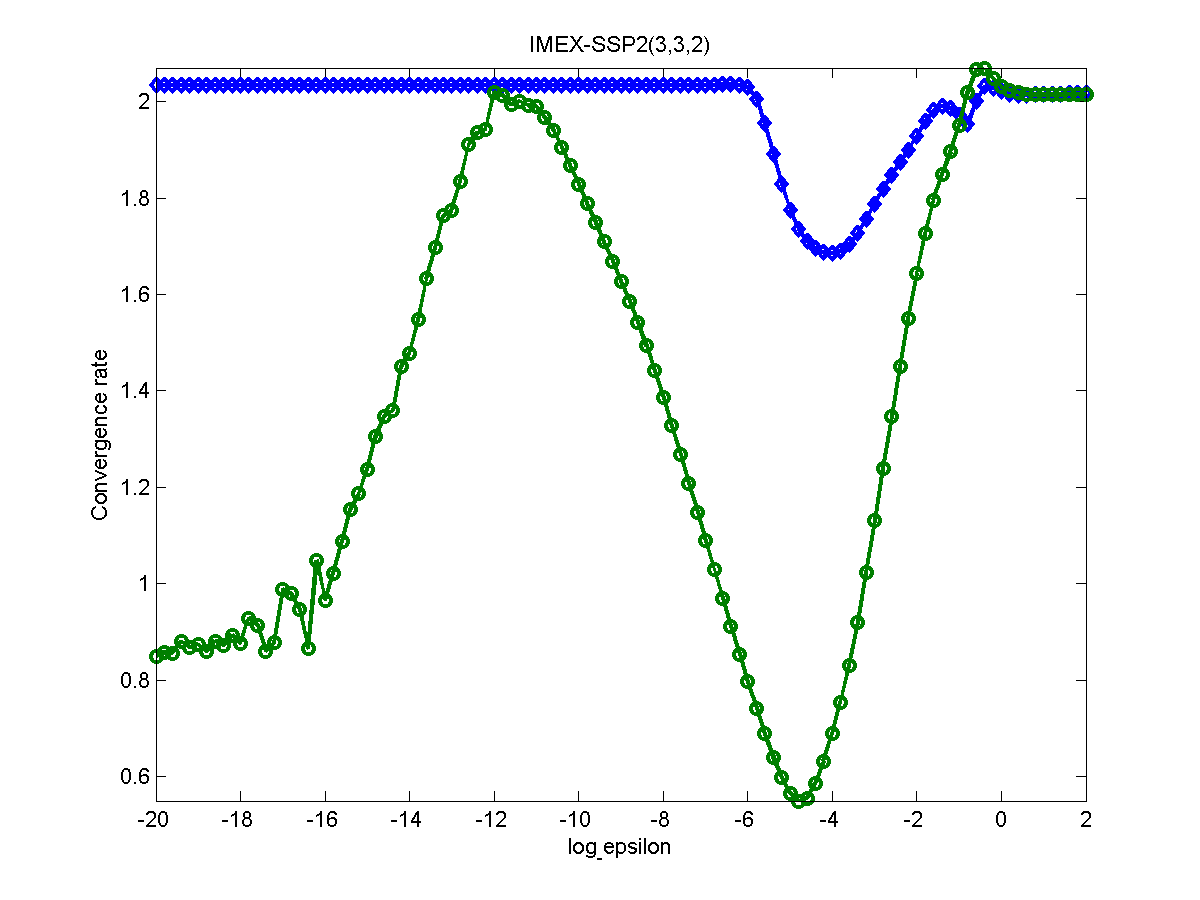}
  \includegraphics[width=0.38\textwidth]{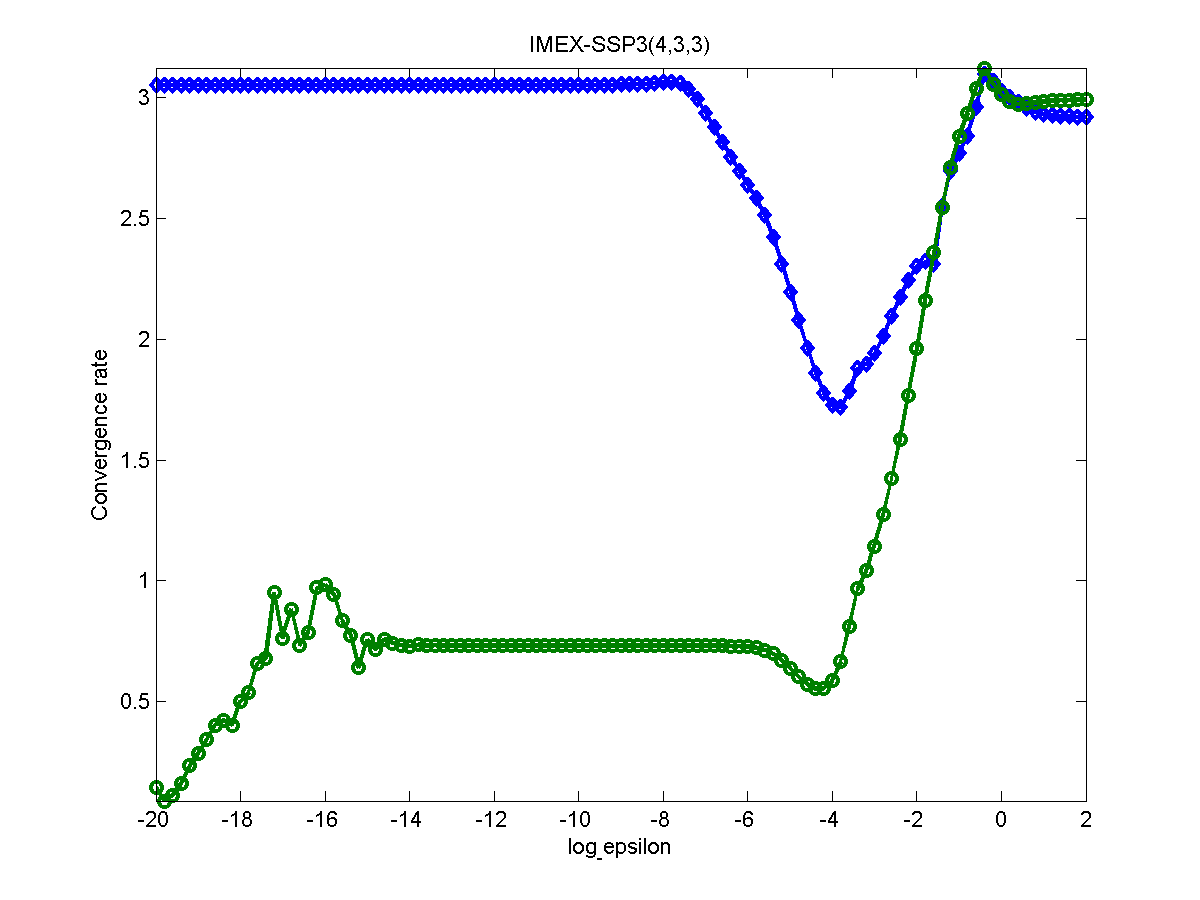}\\
\caption{Convergence behaviors of eight schemes designed in this paper and two from Pareschi and Russo \cite{Pareschi:2005} for Pareschi and Russo's problem with non-equilibrium initial conditions. The blue lines with diamonds are for $x$ and the green lines with circles are for $y$.}\label{fig.ConvergParWtNonEquil}
\end{figure}

Representative multi-dimensional displays for the $L_2$-norm error tested with two selected schemes are presented in Fig. \ref{fig.L2NormError}. The left column shows the $x$ variable's $L_2$ error and the right column shows the $y$ variable's $L_2$ error. The convergence behaviors of eight schemes designed in this paper and two from Pareschi and Russo \cite{Pareschi:2005} are shown in Fig. \ref{fig.ConvergParWtEquil} for equilibrium initial conditions and Fig. \ref{fig.ConvergParWtNonEquil} for non-equilibrium conditions. The convergence rates in Figs. \ref{fig.ConvergParWtEquil} and \ref{fig.ConvergParWtNonEquil} are obtained by fitting $E(\Delta t)$ to a straight line in a double logarithm coordinate. The temporal region for second-order schemes' fitting is $\Delta t\in [10^{-4}, 10^0]$ and $\Delta t\in[10^{-3}, 10^0]$ for third-order schemes. The assumed ``exact'' solution is computed with $\Delta t = 10^{-6}$. There is nearly the same lower limit of the error, order $10^{-11}$, for different schemes. We assume this is caused by the accumulation of the machine truncation error. Indeed, we observed the error increasing when the time step decreased beyond $10^{-6}$.

We first introduce the results for eight schemes designed in this paper. In every colorized contour plot at the bottom of the subplots in Fig. \ref{fig.L2NormError}, there is an overlapping straight slash $\log\Delta t = \log \varepsilon$. The location of this slash indicates where the ``ridge'' of the error is located. Because of this ridge, the convergence rate in Figs. \ref{fig.ConvergParWtEquil} and \ref{fig.ConvergParWtNonEquil} has a ``dip'' in the intermediate region of $\varepsilon$. The fact accompanying this is that the matching degree of the straight line is very low; hence, the convergence rate is poorly defined in this region. However, the convergence rate in this region is helpful for semi-quantitatively observing the deterioration of convergence in the intermediate stiffness, which is a common phenomenon observed for many IMEX schemes \cite{Kennedy:2003, Pareschi:2001} and appears when $\log\Delta t\simeq \log\varepsilon$, that is $\Delta t=O(\varepsilon)$, as demonstrated by the reference slash line at the bottom of the subplots in Fig. \ref{fig.L2NormError}. A common feature of this deterioration is that it is more severe for the second variable $y$ than the first variable $x$ \cite{Kennedy:2003, Pareschi:2001}, though of a different degree for different initial conditions and schemes. Despite the reduction of the order of accuracy that appears in the intermediate region, complete recovery appears on two sides for the equilibrium initial conditions in Fig. \ref{fig.ConvergParWtEquil}. The leftmost point of the $\varepsilon$ axis is in fact zero; these schemes behave in the same manner in infinite relaxation as they do in finite relaxation. For non-equilibrium cases in Fig. \ref{fig.ConvergParWtNonEquil}, however, only the first (4,3,2), third (4,3',2), and sixth (6,4,3) schemes can completely recover to the designed accuracy order on the left-hand side of the $\varepsilon$ axis, whereas the other five merely recover convergence accuracy to first-order. These observations seem to indicate that the first (4,3,2), third (4,3',2), and sixth (6,4,3) schemes precede the others. However, the second (3',3,2) and fifth (4',4,2) schemes perform much better for equilibrium cases; the ``ridge'' becomes a negligible ``ripple'', as illustrated in Fig. \ref{fig.L2NormError}. A common feature of these two schemes is the monotony of the internal sub-time level; however, we do not have conclusive evidence that this is linked to the good convergence behavior.

Now, we briefly introduce the results for two schemes from Pareschi and Russo \cite{Pareschi:2005} as a comparison. IMEX-SSP2(3,3,2) is the best second-order scheme and IMEX-SSP3(4,3,3) is the only third-order scheme from Pareschi and Russo \cite{Pareschi:2005}. Both schemes can completely recover to the designed accuracy order on the left-hand side of the $\varepsilon$ axis for the first variable $x$, but cannot for the second variable $y$ regardless of the initial condition type. In particular, when $\varepsilon = 0$, both schemes cannot be used (the leftmost point of the ¦Å axis is not zero for these two schemes). As a comparison, our new schemes are indeed capable in the zero relaxation limit without degradation of convergence, as demonstrated by our tests. If a computer code based on, for example, relaxation magnetohydrodynamics (MHD) and integrated with an IMEX-SSP-like scheme is aimed at running from a relaxation range to an ideal range ($\varepsilon = 0$, corresponding to ideal MHD), the source code may be written as an if-else statement for the case of $\varepsilon=0$. However, if it is integrated with our new scheme, the if-else statement is unnecessary. Thus, our new schemes can be used in the zero relaxation limit in a unified and convenient manner. In our experience, the source code would be between half and a quarter of the length of the code that includes if-else statements.

\begin{figure}[htbp]
\centering
  \includegraphics[width=0.38\textwidth]{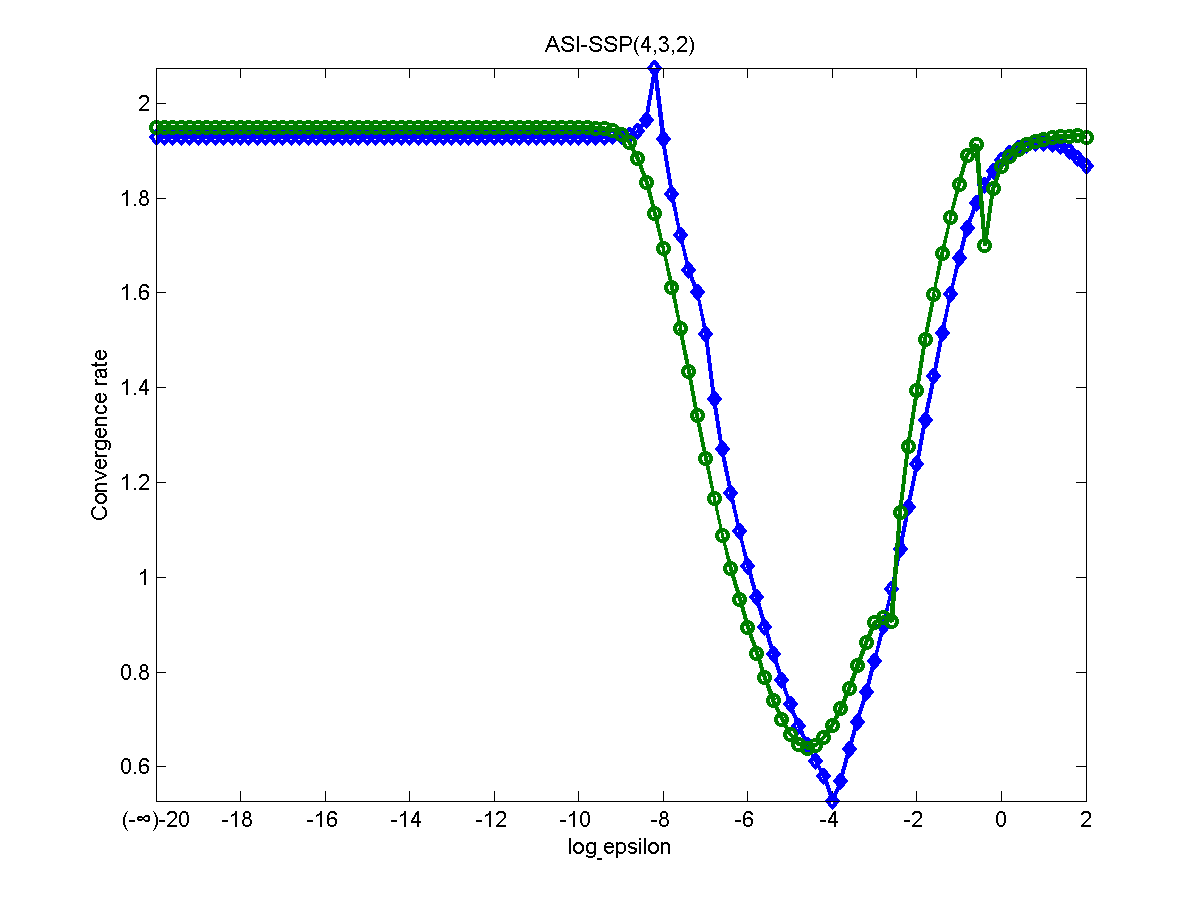}
  \includegraphics[width=0.38\textwidth]{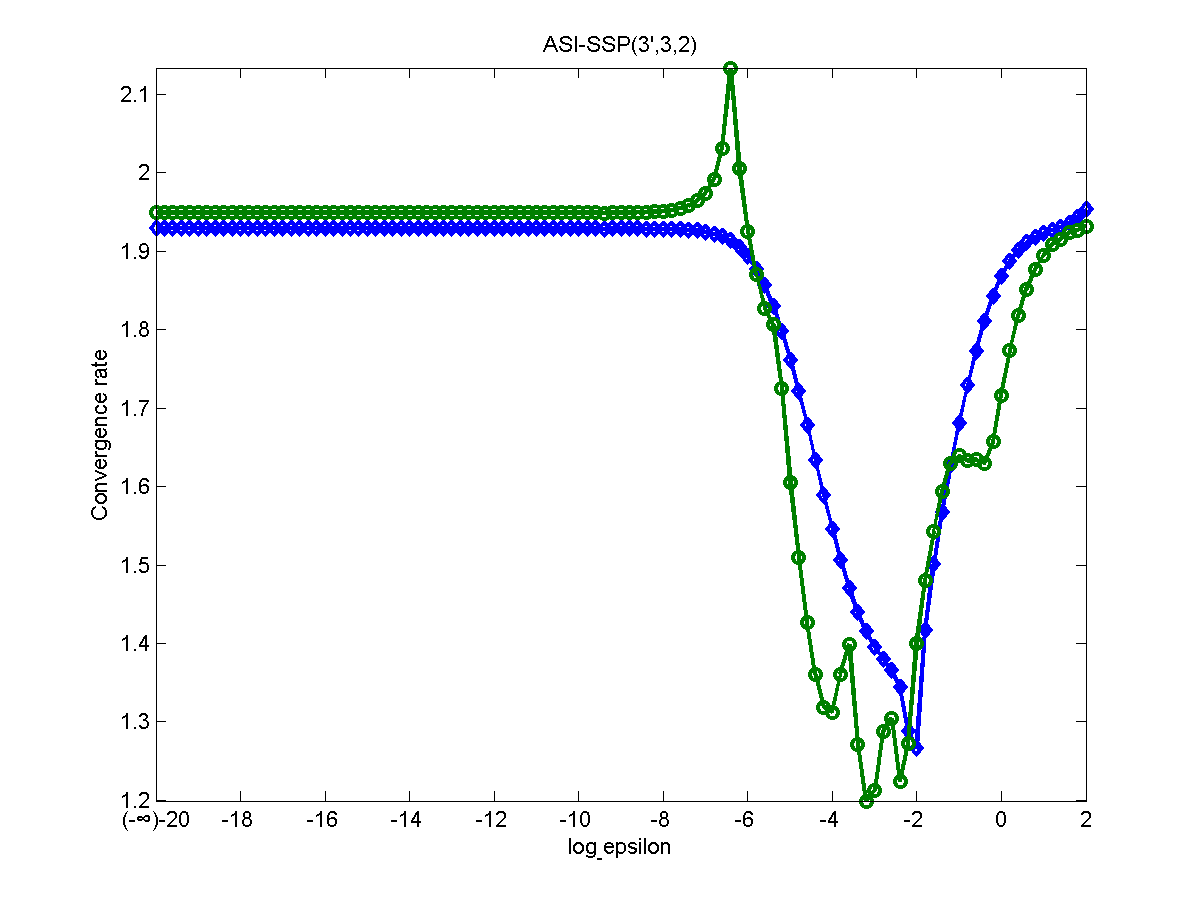}\\
  \includegraphics[width=0.38\textwidth]{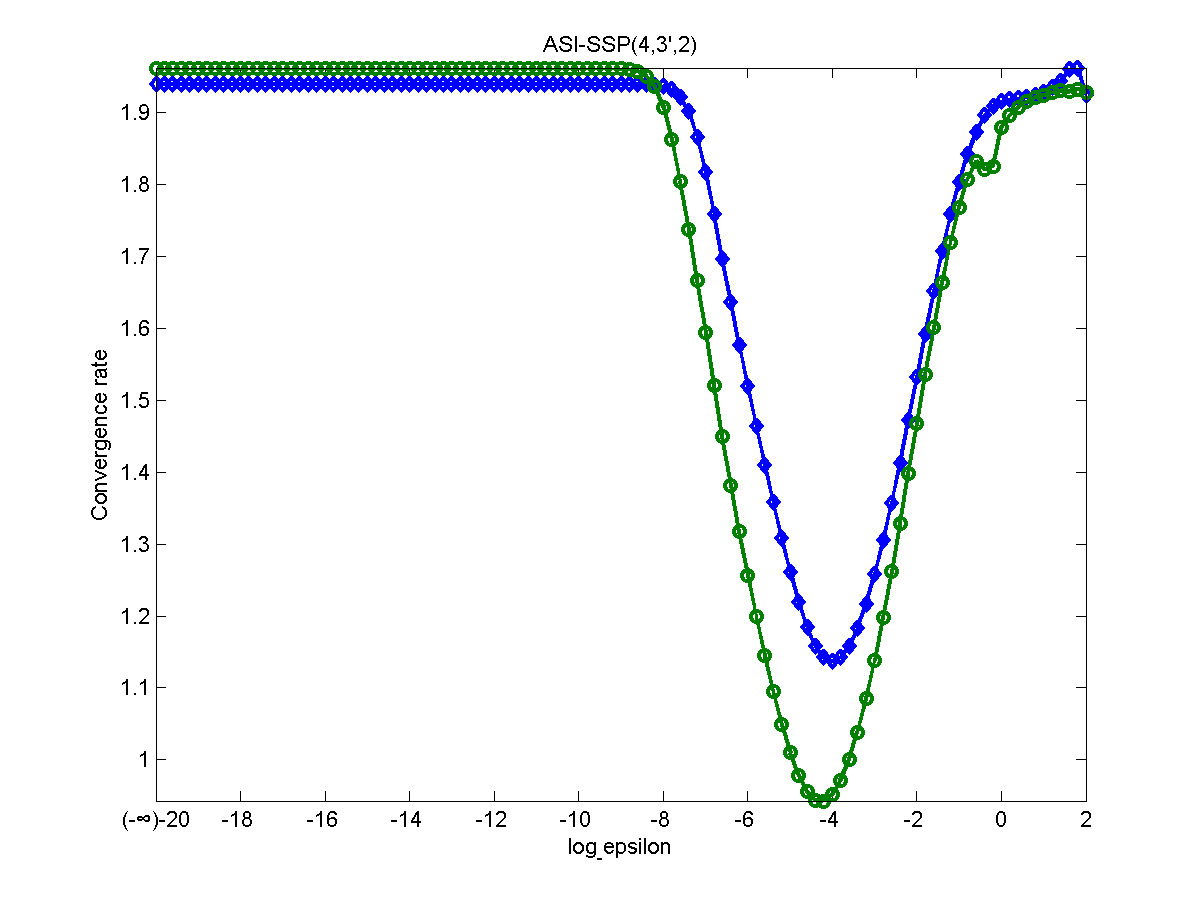}
  \includegraphics[width=0.38\textwidth]{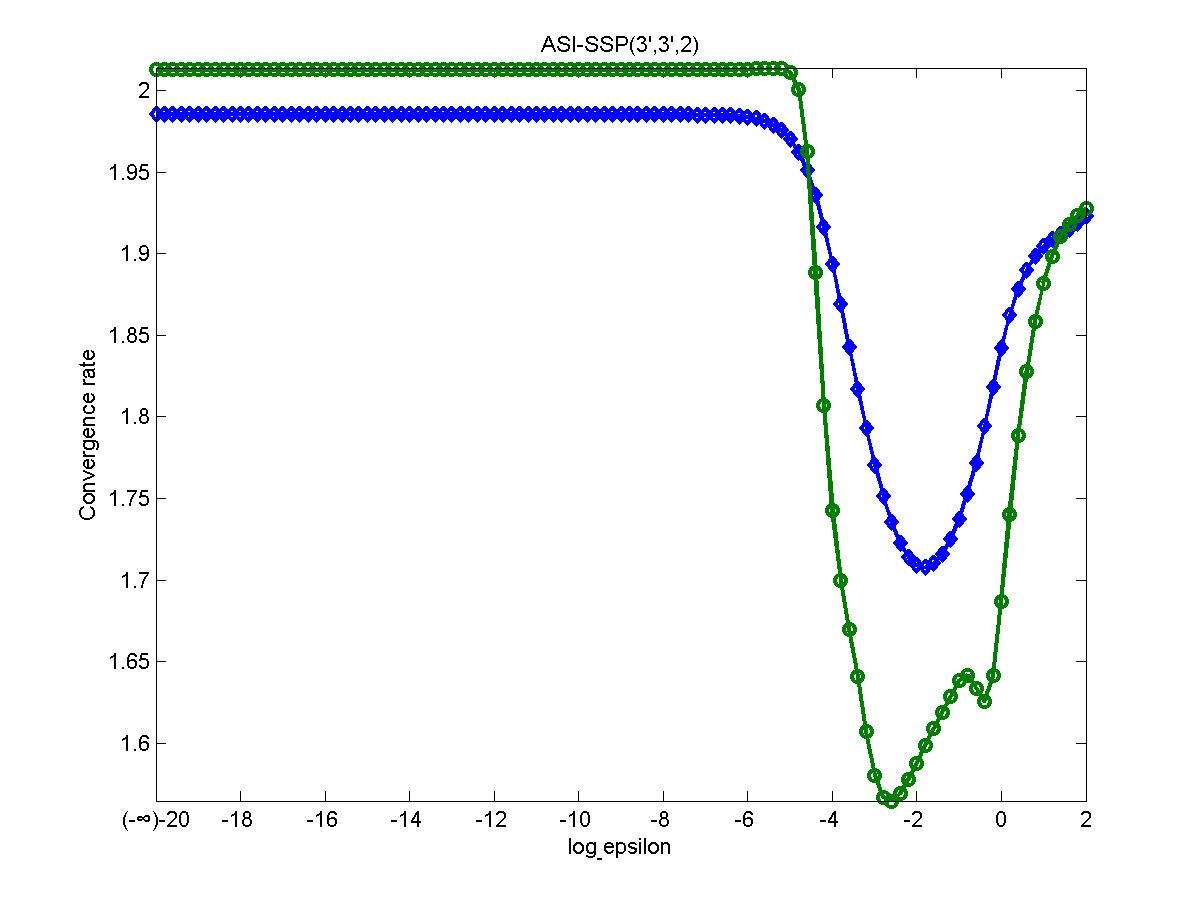}\\
  \includegraphics[width=0.38\textwidth]{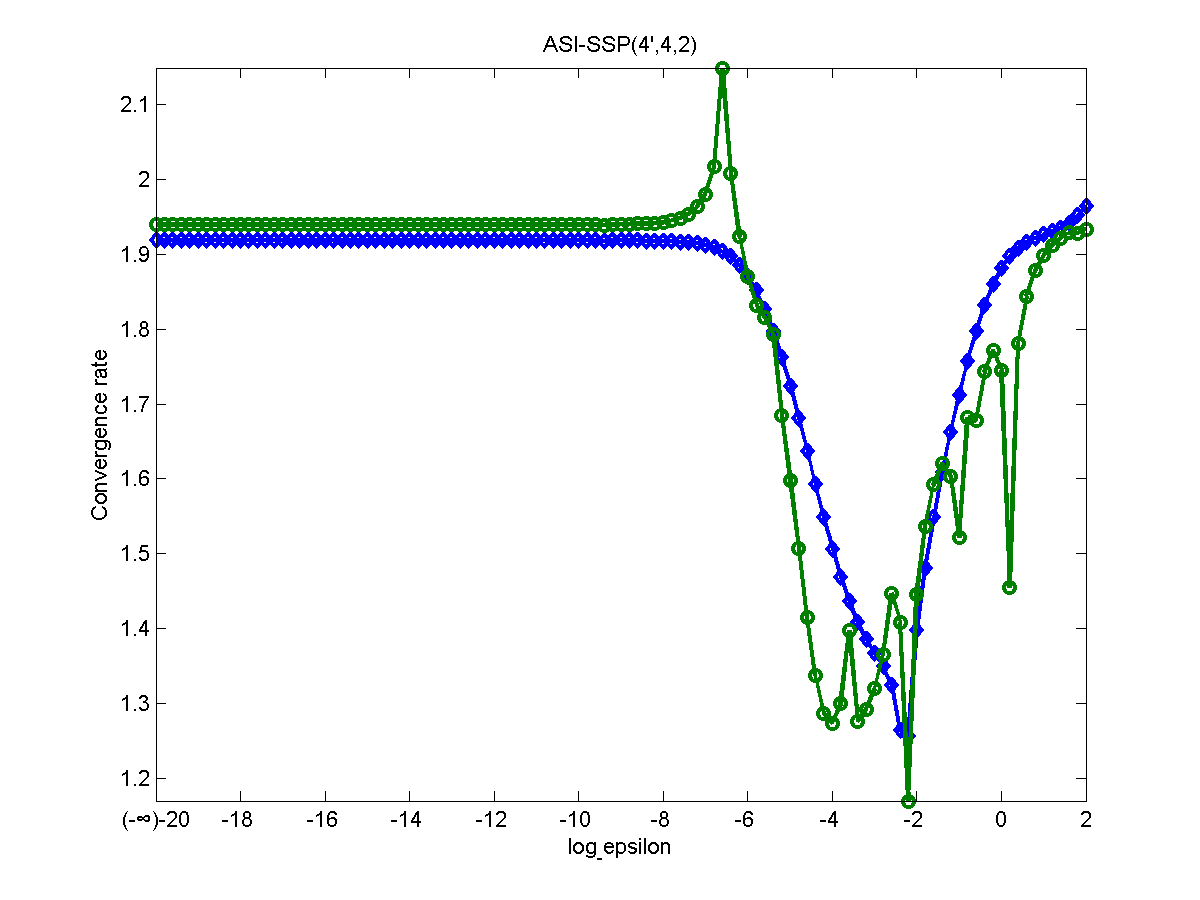}
  \includegraphics[width=0.38\textwidth]{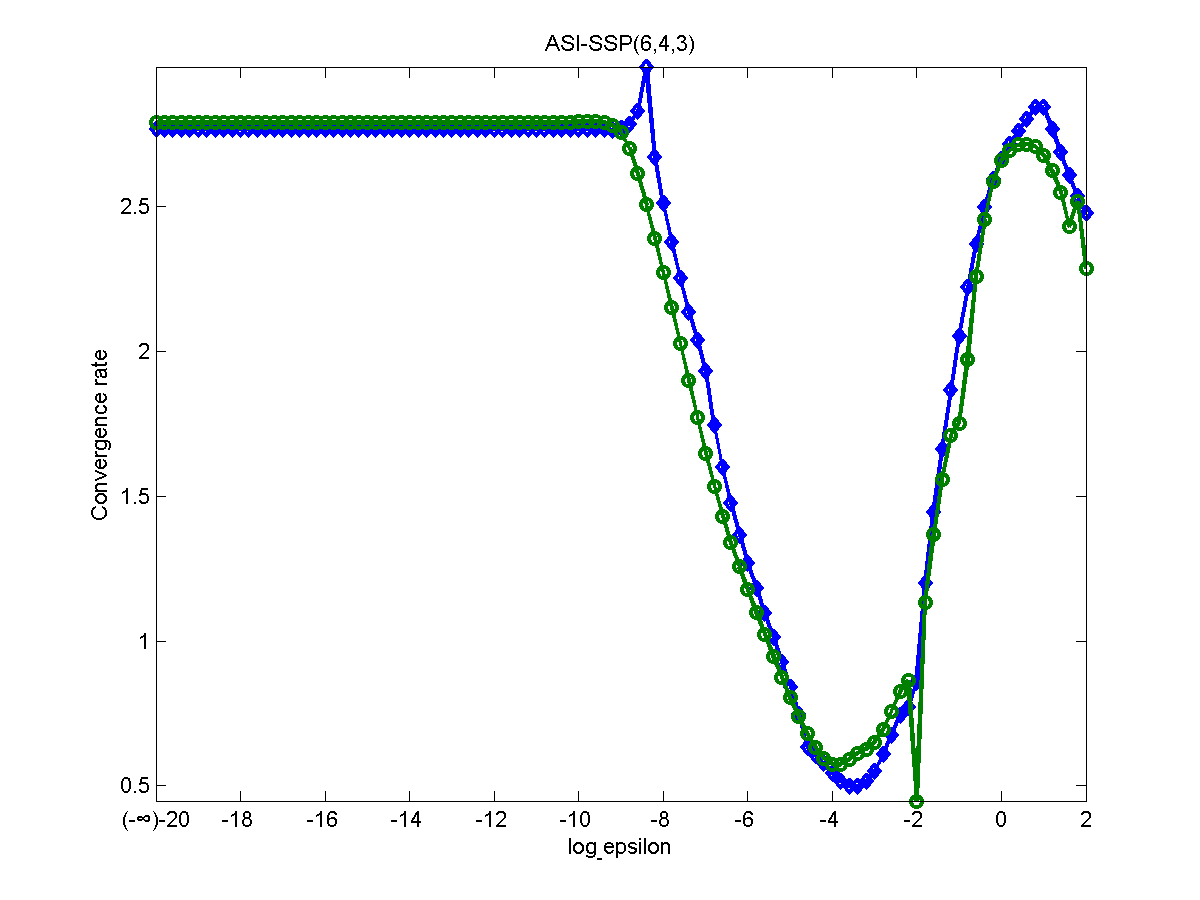}\\
  \includegraphics[width=0.38\textwidth]{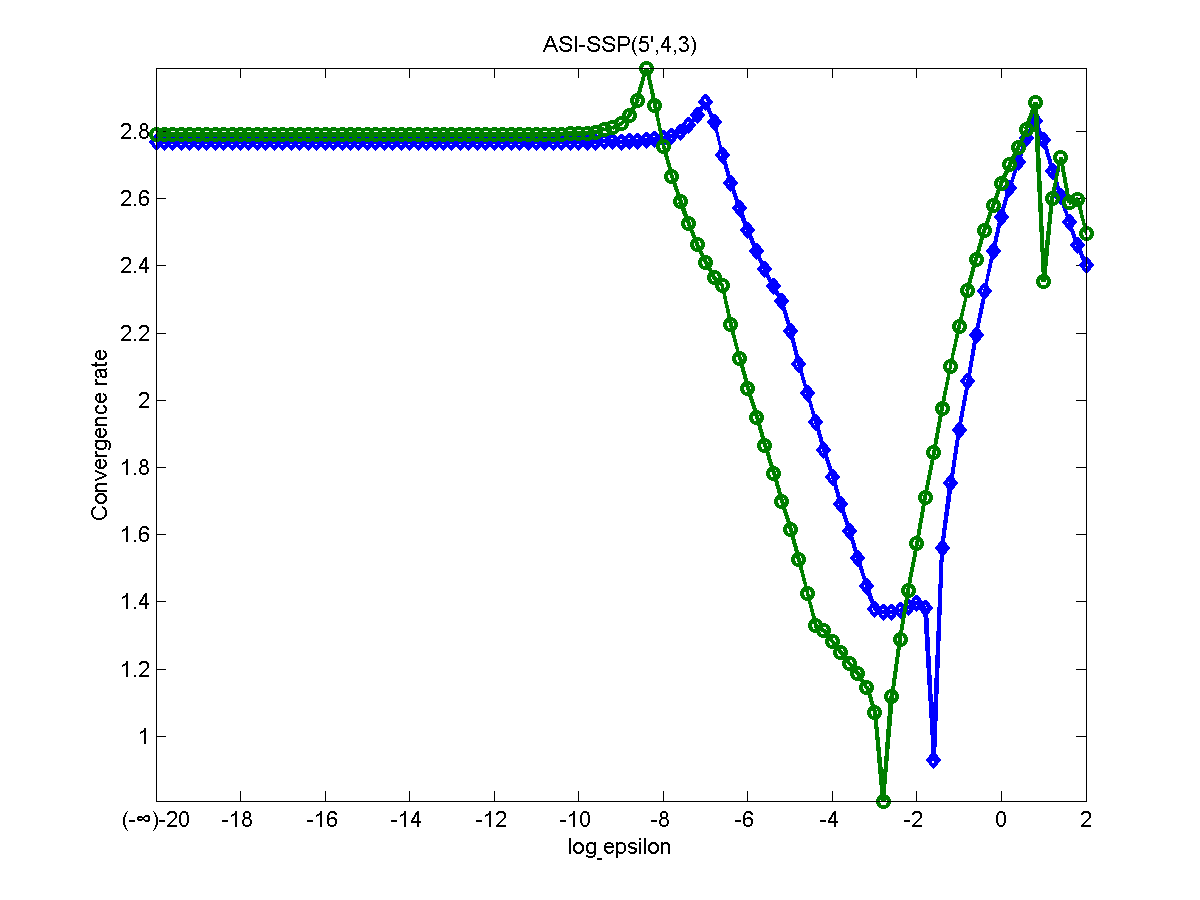}
  \includegraphics[width=0.38\textwidth]{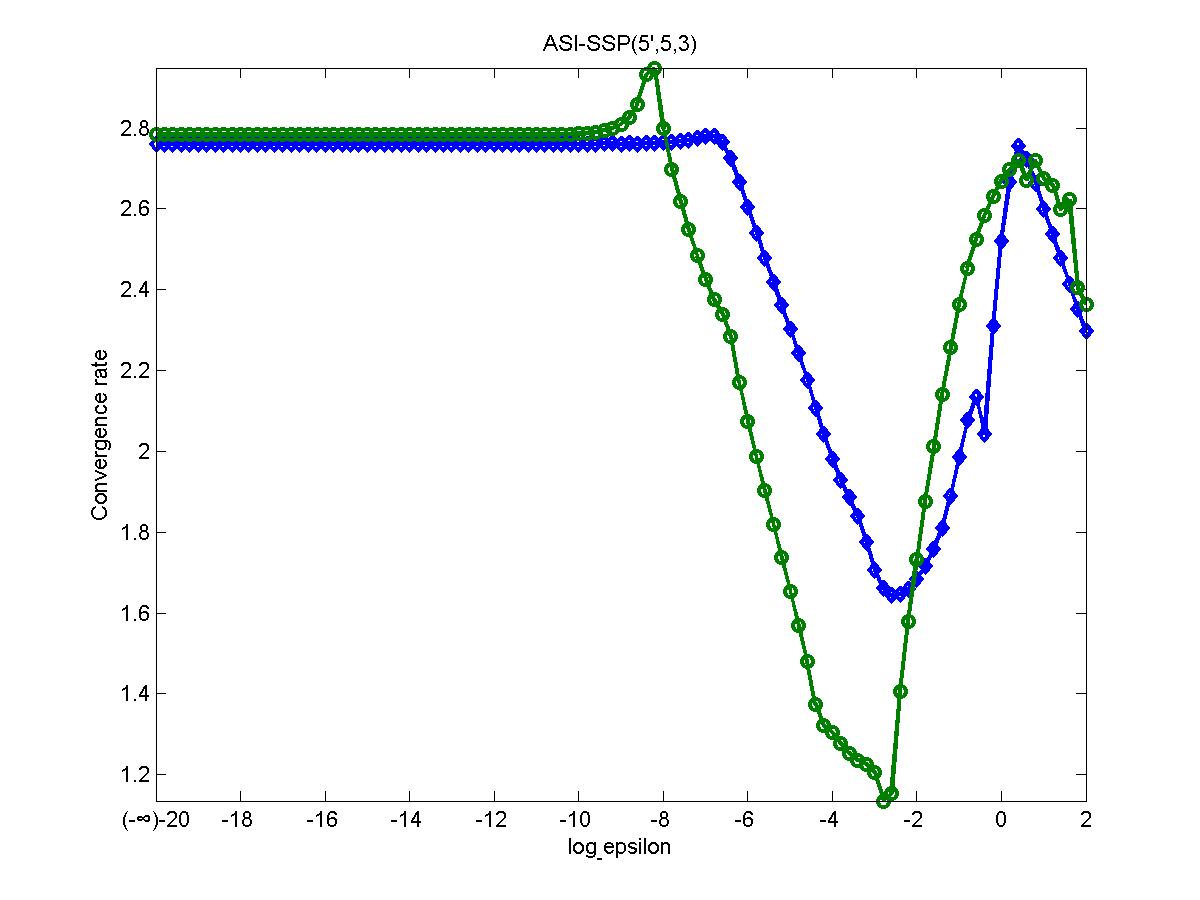}\\
  \includegraphics[width=0.38\textwidth]{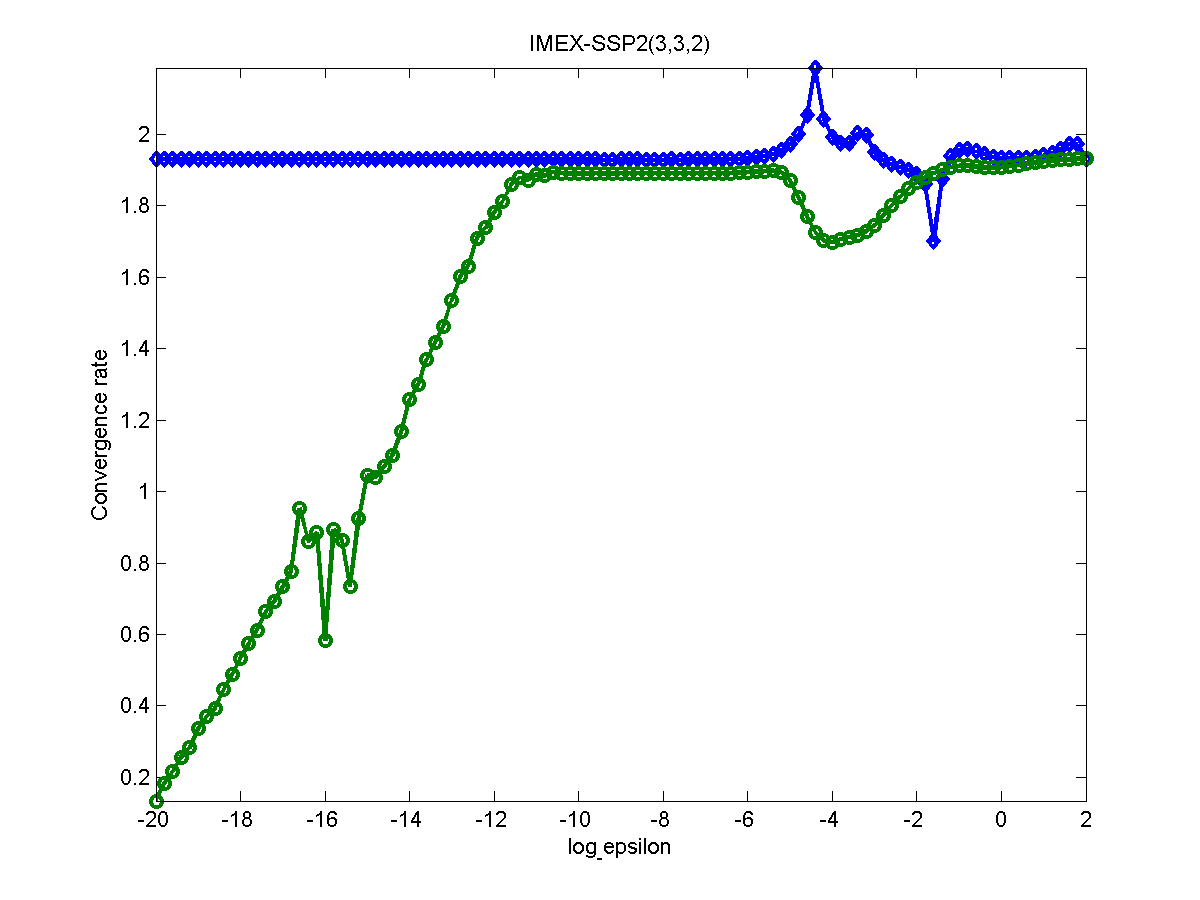}
  \includegraphics[width=0.38\textwidth]{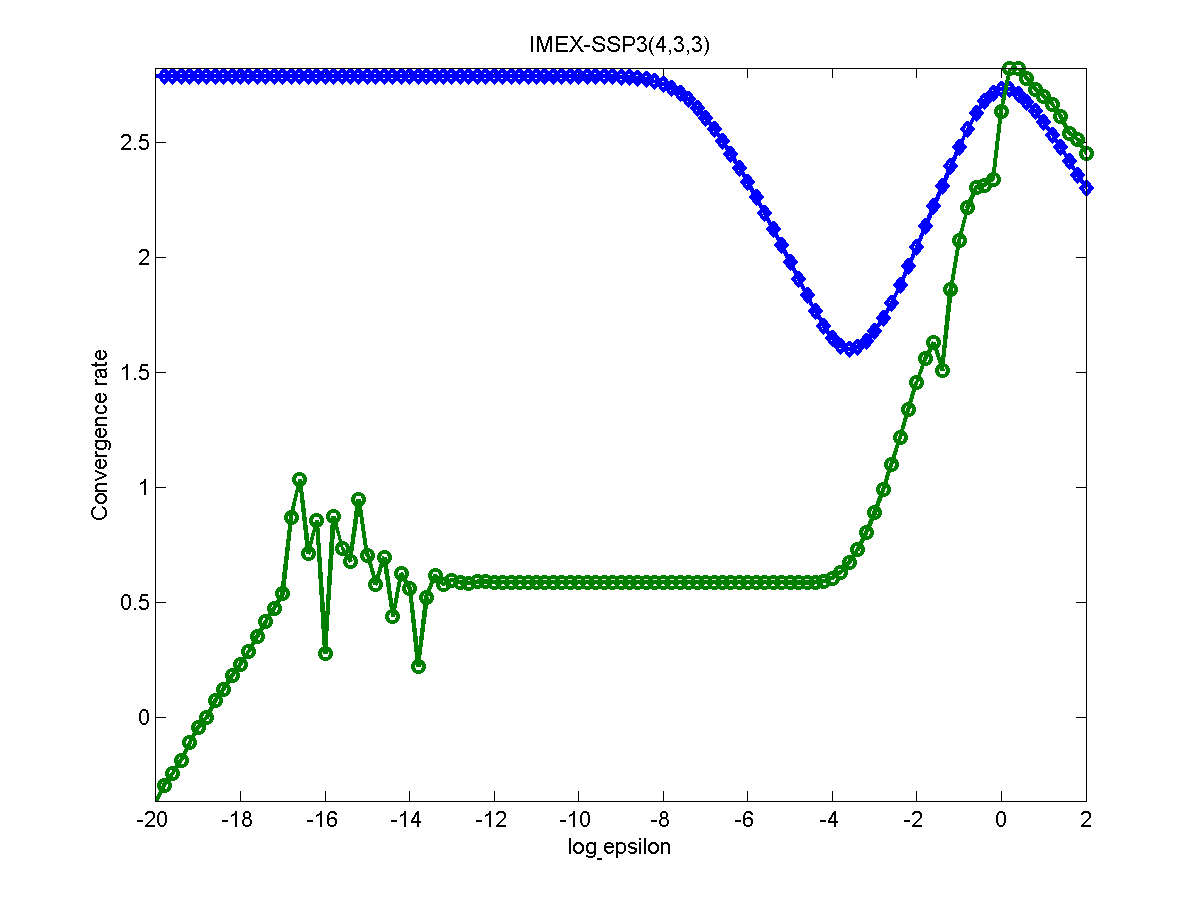}\\
\caption{Convergence behaviors of eight schemes designed in this paper and two from Pareschi and Russo \cite{Pareschi:2005} for van der Pol's equation with equilibrium initial conditions. The blue lines with diamonds are for $x$ and the green lines with circles are for $y$.}\label{fig.ConvergVanWtEquil}
\end{figure}

\begin{figure}[htbp]
\centering
  \includegraphics[width=0.38\textwidth]{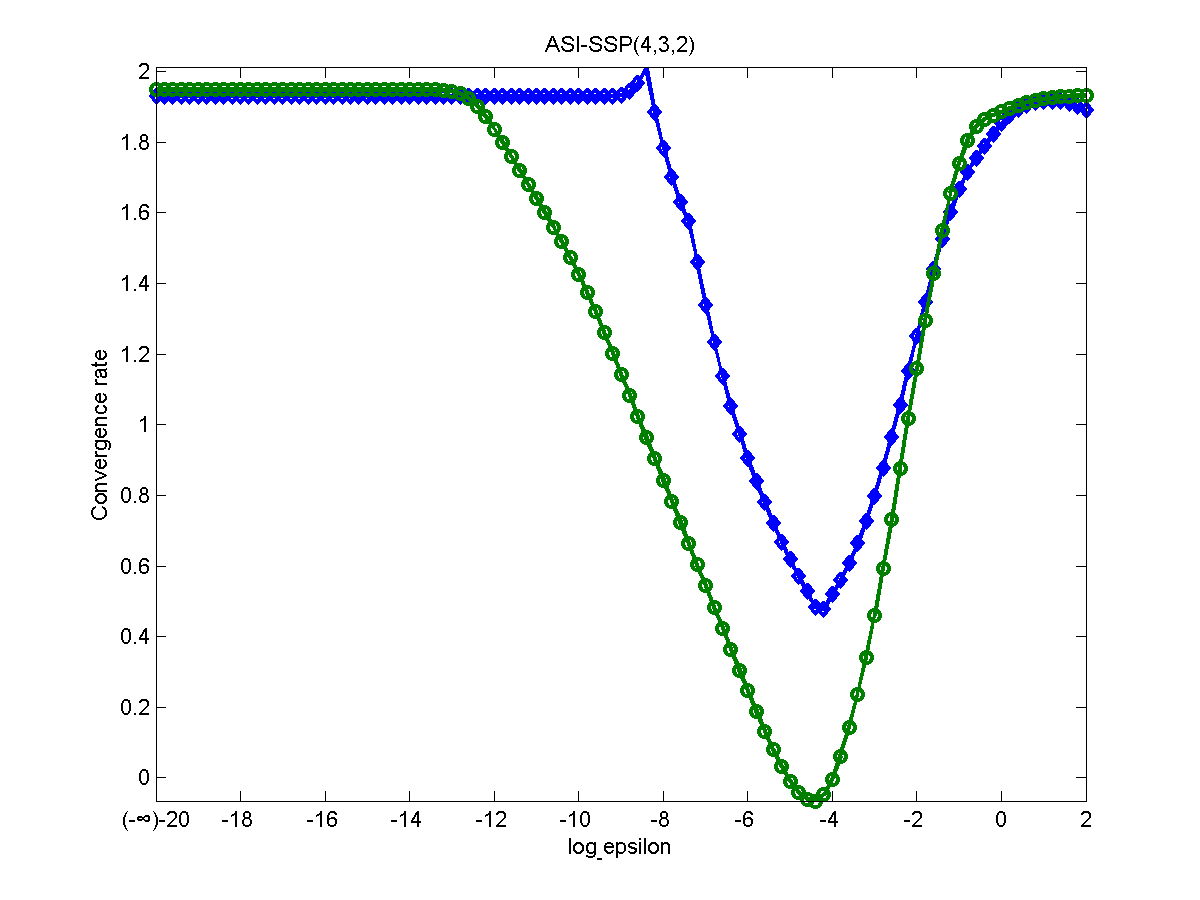}
  \includegraphics[width=0.38\textwidth]{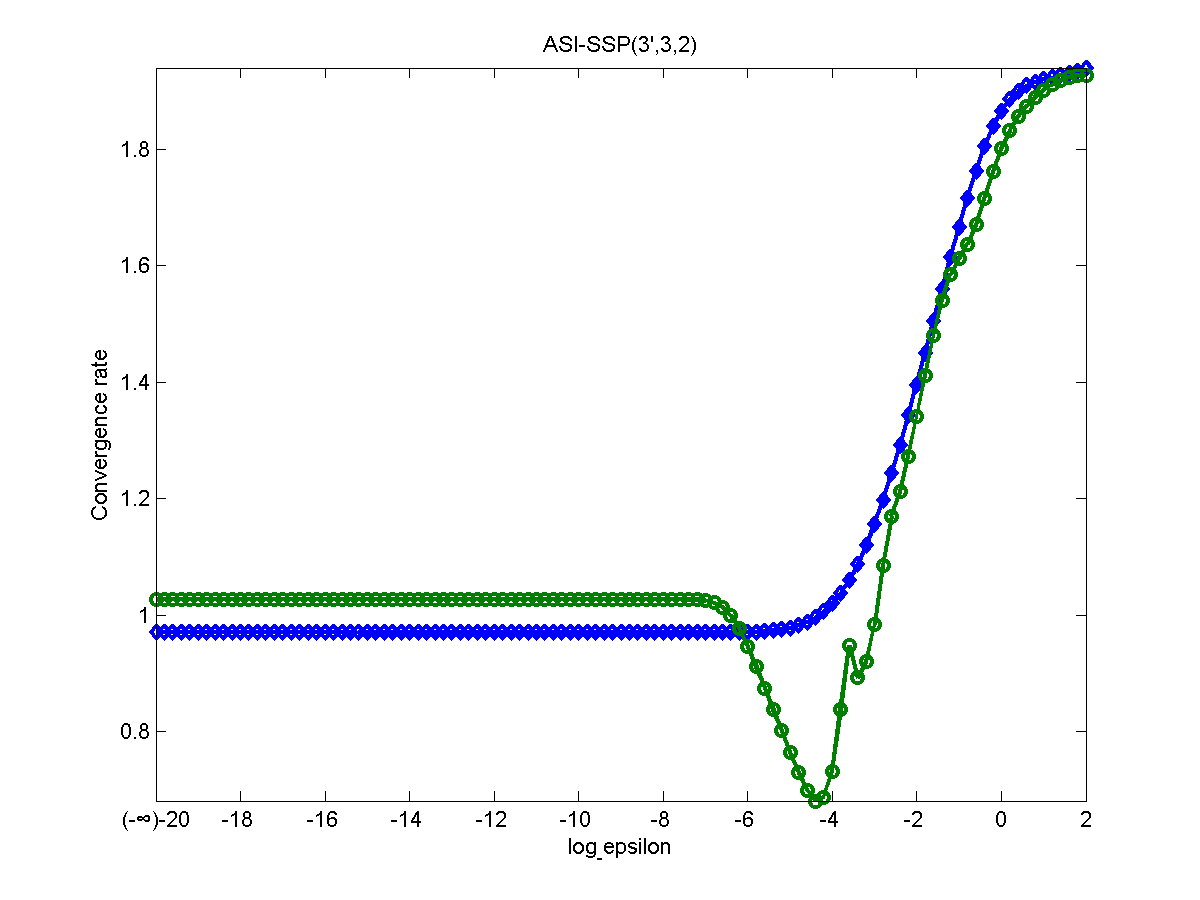}\\
  \includegraphics[width=0.38\textwidth]{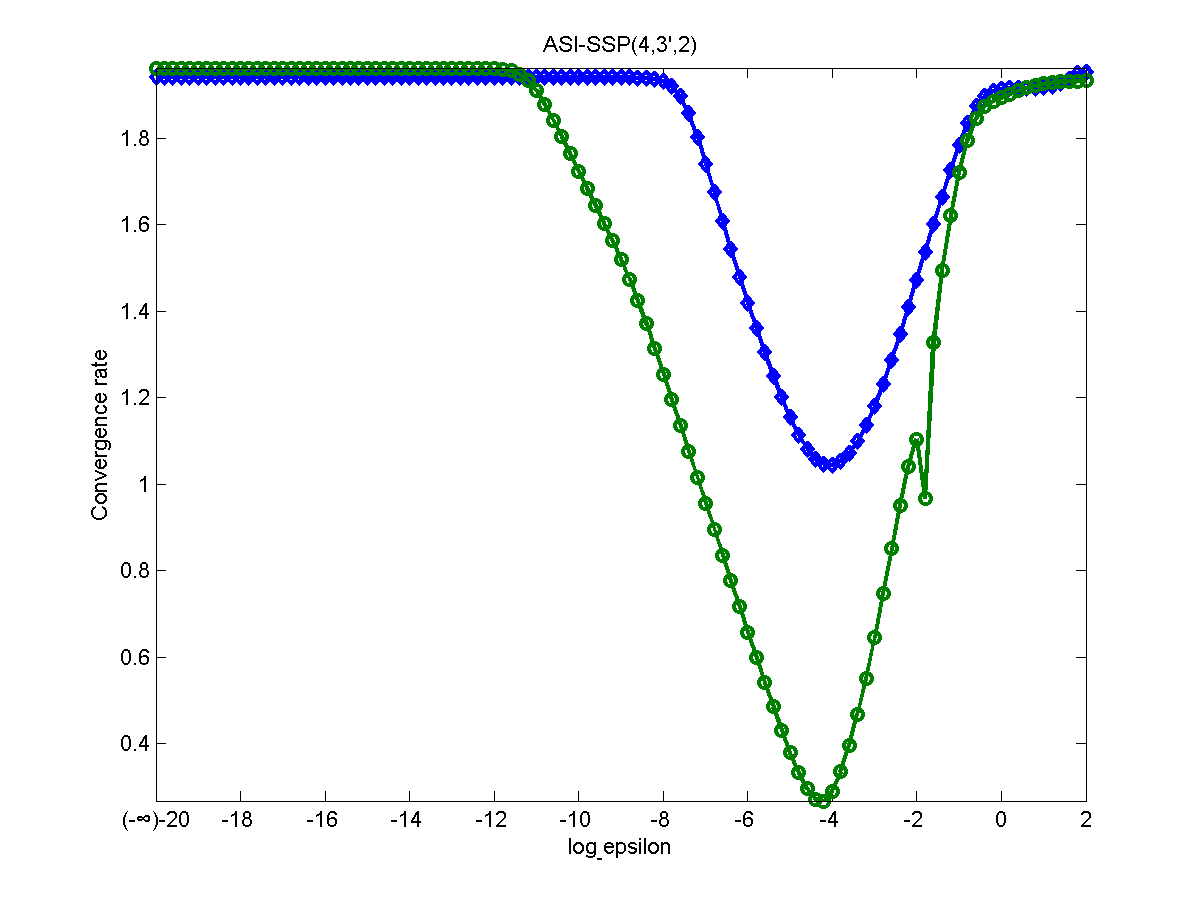}
  \includegraphics[width=0.38\textwidth]{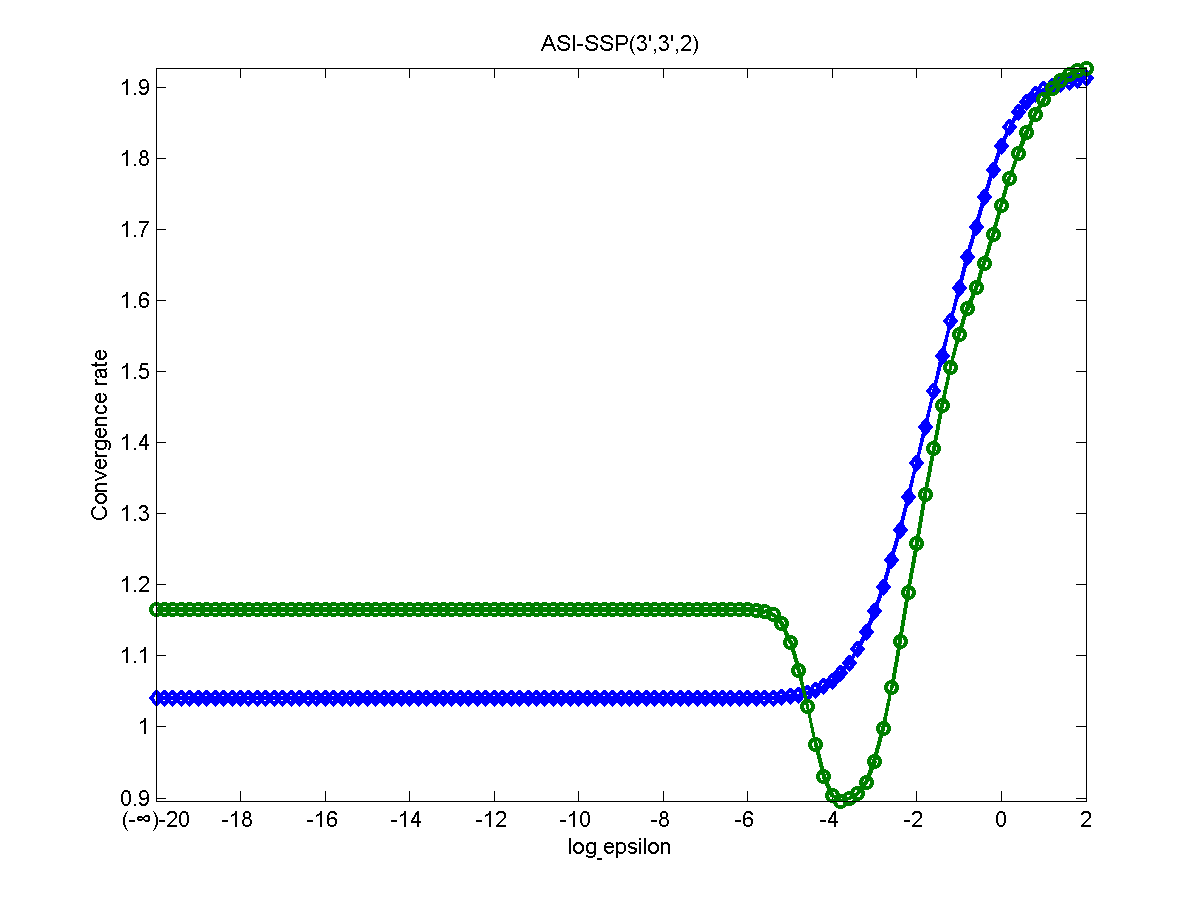}\\
  \includegraphics[width=0.38\textwidth]{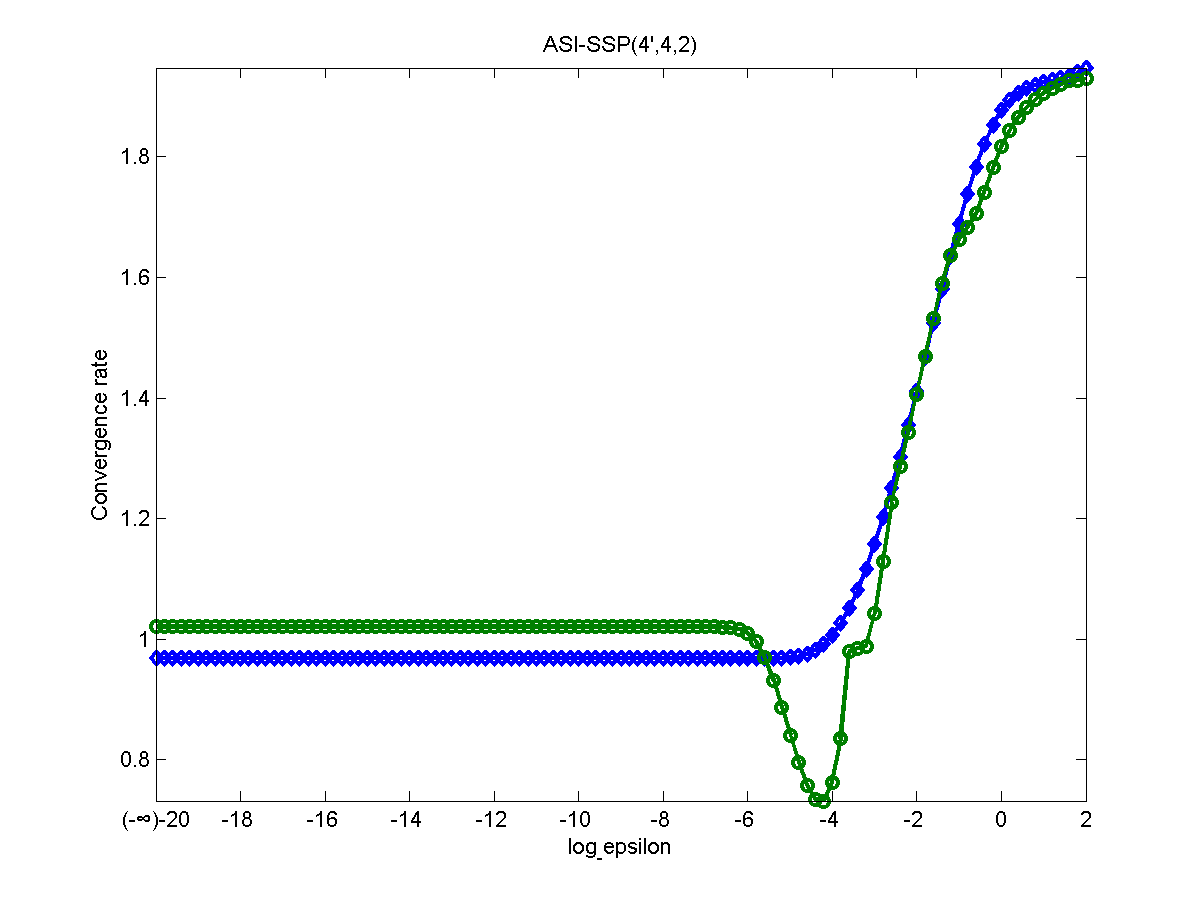}
  \includegraphics[width=0.38\textwidth]{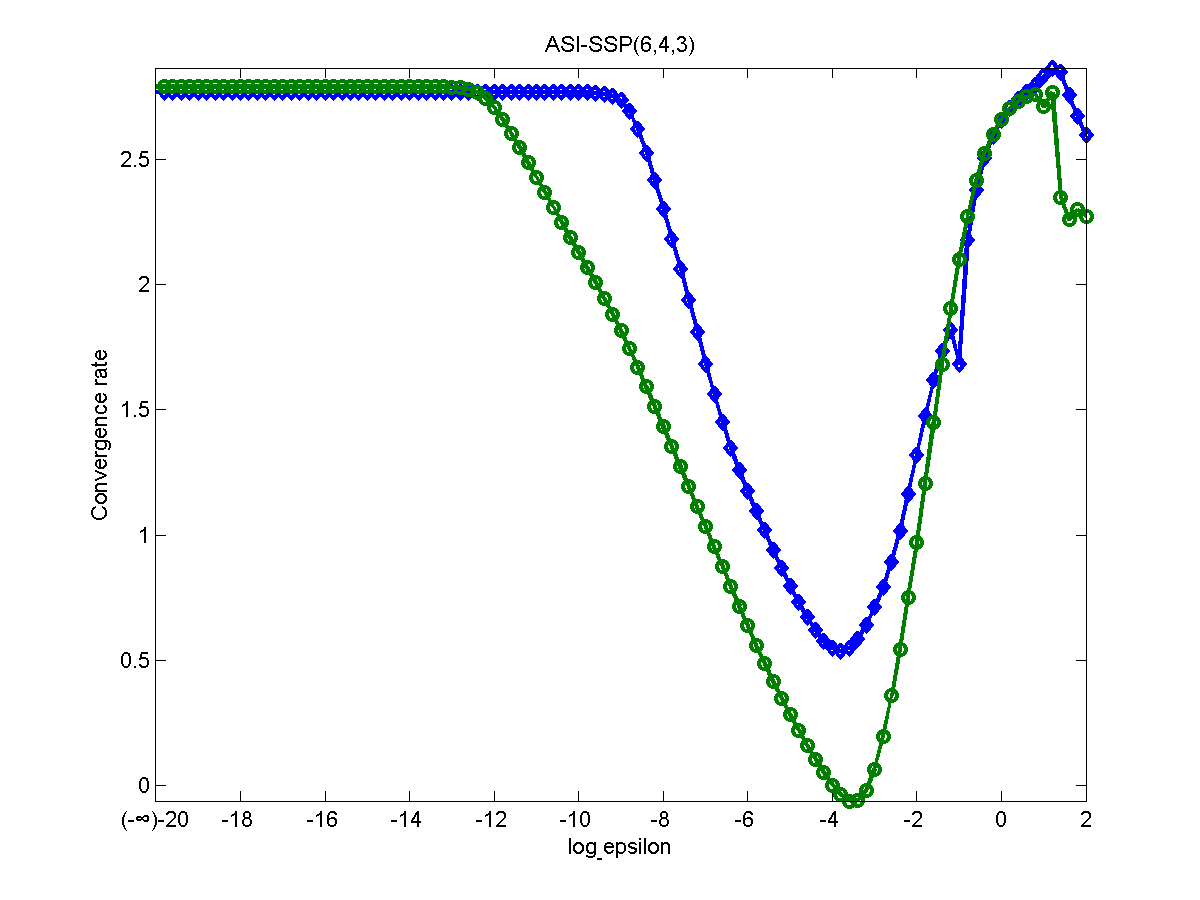}\\
  \includegraphics[width=0.38\textwidth]{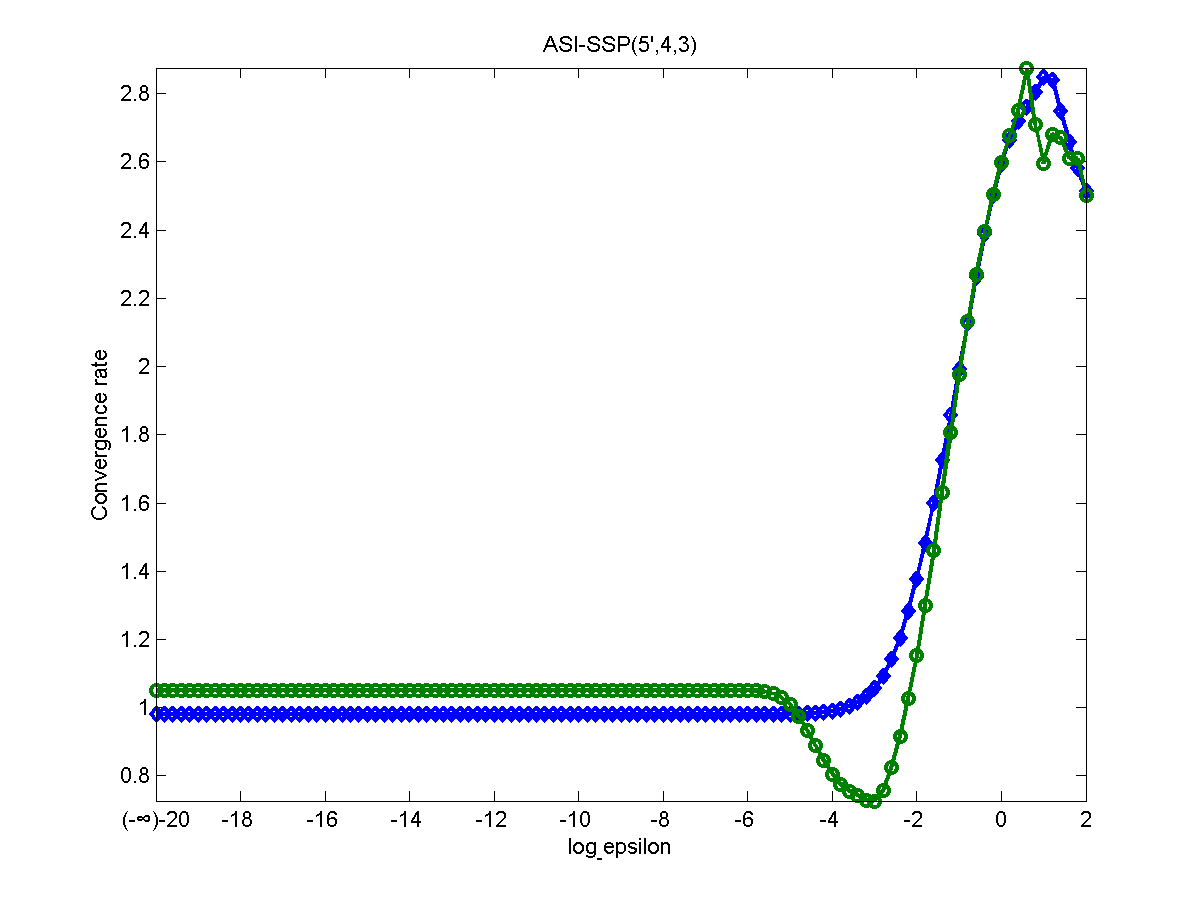}
  \includegraphics[width=0.38\textwidth]{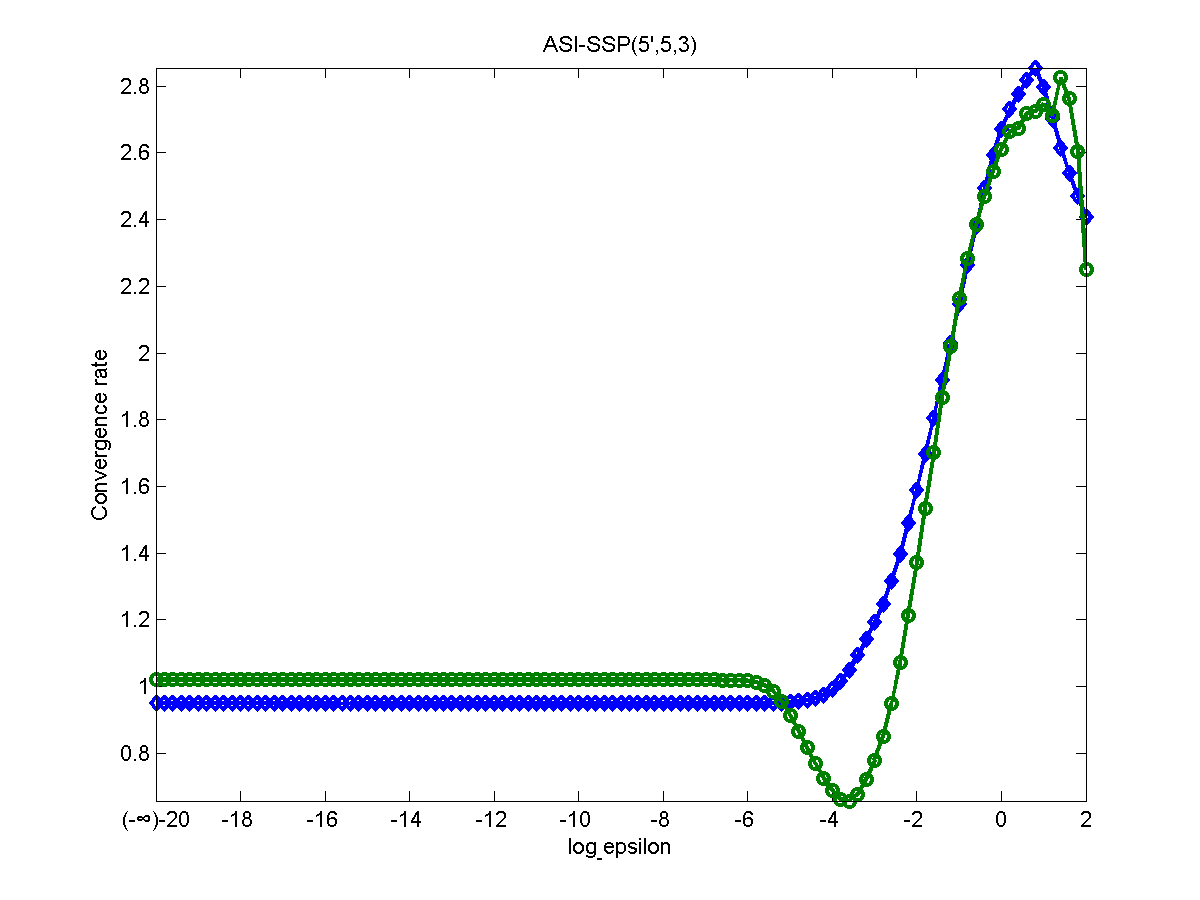}\\
  \includegraphics[width=0.38\textwidth]{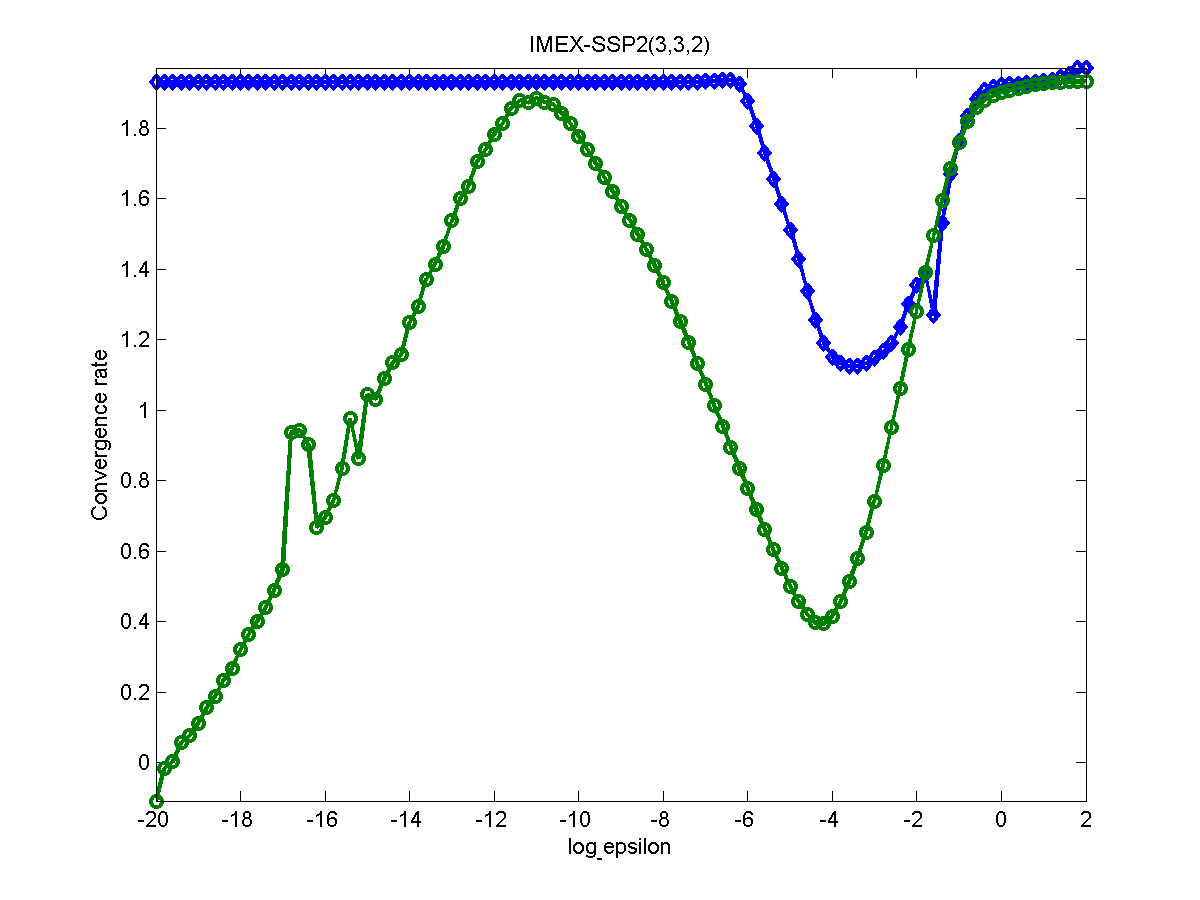}
  \includegraphics[width=0.38\textwidth]{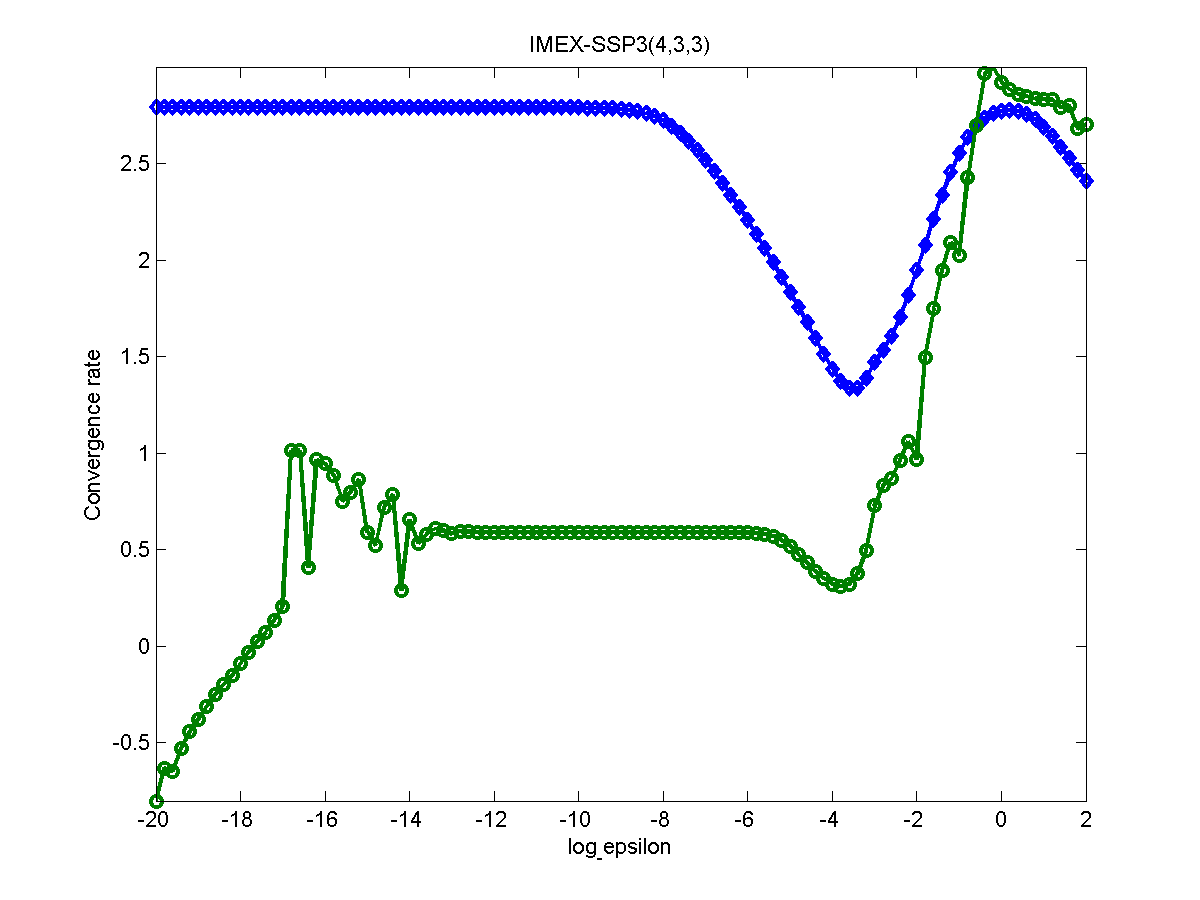}\\
\caption{Convergence behaviors of eight schemes designed in this paper and two from Pareschi and Russo \cite{Pareschi:2005} for van der Pol's equation with non-equilibrium initial conditions. The blue lines with diamonds are for $x$ and the green lines with circles are for $y$.}\label{fig.ConvergVanWtNonEquil}
\end{figure}

Finally, we consider van der Pol's equation \cite{Kennedy:2003} to test the convergence behavior:
\begin{equation}
\dot{x}(t) = y(t), \dot{y}(t) = ((1-x(t)^2)y(t) - x(t))/\varepsilon.
\end{equation}
As before, the terms divided by $\varepsilon$ are integrated with the implicit method, whereas the other terms are integrated explicitly. Initial conditions are considered in two forms: equilibrium initial conditions accomplished with $x(0)=2, y(0)=-2/3,$ and non-equilibrium conditions specified by replacing the condition on $y$ with $y(0)=-1$. The integration time interval is $t\in[0,5]$. The other processes are as stated previously. The results presented in Figs. \ref{fig.ConvergVanWtEquil} and \ref{fig.ConvergVanWtNonEquil} are very similar to those in 
Figs. \ref{fig.ConvergParWtEquil} and \ref{fig.ConvergParWtNonEquil}, thus are not repeated here.

\section{Conclusion}\label{sec.conclusion}

We construct eight IMEX RK schemes up to third order of the type in which all stages are implicit so that they can be used in the zero relaxation limit in a unified and convenient manner. All schemes attain the SSP property in the limiting case in the same manner as Pareschi and Russo's IMEX-SSP schemes, and two schemes, the first (4,3,2) and third (4,3',2), are SSP not only for the explicit part but also the implicit part and the entire IMEX scheme. The third (4,3',2) and sixth (6,4,3) schemes have a property that the stable region contains an interval on the imaginary axis. The first (4,3,2), third (4,3',2), and sixth (6,4,3) schemes can completely recover to the designed accuracy order in two sides of the relaxation parameter for both equilibrium and non-equilibrium initial conditions, although the reduction of the accuracy order appears in the intermediate region. However, the second (3',3,2) and fifth (4',4,2) schemes achieve much better convergence for equilibrium cases.

{\em Acknowledgements:} This work is in part supported by NSFC Grant Nos. 11405167, 51407171, 11571293 and 11505171, Project 2015B0201023 of CAEP, and Key Laboratory of Pulsed Power, CAEP, Contract No. PPLF2014PZ07.




\bibliographystyle{elsarticle-num}
\bibliography{<your-bib-database>}



\end{document}